\documentclass[11pt,a4paper,reqno]{amsart}
\usepackage{times} 
\usepackage{amsmath} 
\usepackage{amssymb}  
\usepackage{amsthm}
\usepackage{latexsym}
\usepackage{amsfonts,bbm}
\usepackage{xcolor}
\usepackage{mathtools}
\usepackage{enumerate}
\usepackage{cite}
\usepackage{tikz}
\usepackage{graphicx,caption}
\usepackage{todonotes}
\usepackage{float}
\usepackage{nicefrac}
\usepackage{hyperref}
\usepackage{xcolor}
\hypersetup{
    colorlinks,
    allcolors = {green!50!black},
    urlcolor={red!50!black}
}
\usepackage[nameinlink]{cleveref}
\crefformat{equation}{\textup{#2(#1)#3}}

\newtheorem{theorem}{Theorem}[section] 

\newtheorem{remark}[theorem]{Remark}
\newtheorem{proposition}[theorem]{Proposition} 
\newtheorem{corollary}[theorem]{Corollary}  
\newtheorem{lemma}[theorem]{Lemma} 
\newtheorem{defi}[theorem]{Definition}

\crefname{hypothesis}{Hypothesis}{Hypotheses}
\crefname{subsection}{Subsection}{Subsections}
\crefname{section}{Section}{Sections}
\crefname{theorem}{Theorem}{Theorems}
\crefname{lemma}{Lemma}{Lemmas}
\crefname{remark}{Remark}{Remark}
\crefname{appendix}{Appendix}{Appendices}
\crefname{proposition}{Proposition}{Propositions}
\crefname{table}{Table}{Tables}
\crefname{figure}{Figure}{Figures}

\usepackage[foot]{amsaddr}
 
\newcommand{\R}{\mathbb R}
\newcommand{\domY}{{Y}}
\newcommand{\domX}{{X}}
\newcommand{\native}{\mathcal N}
\newcommand{\floor}[1]{\left\lfloor #1 \right\rfloor}
\newcommand{\norm}[1]{\left\lVert #1 \right \rVert}
\newcommand{\abs}[1]{\left\lvert #1 \right \rvert}
\newcommand{\inv}[1]{{#1}^{-1}}

\newcommand{\ceil}[1]{\lceil #1 \rceil}

\usepackage{amsopn}

\DeclareMathOperator{\spn}{span}
\DeclareMathOperator{\koopman}{\mathcal K_A}

\newcommand{\paragraf}{\textsection}

\newcommand{\one}{\mathbbm{1}}
\newcommand{\braces}[1]{{\rm (}#1{\rm )}}

\newcommand{\<}{\langle}
\renewcommand{\>}{\rangle}

\newcommand{\wt}{\widetilde}
\newcommand{\wh}{\widehat}
\newcommand{\dist}{\operatorname{dist}}

\newcommand{\N}{\ensuremath{\mathbb N}}    
\newcommand{\Z}{\ensuremath{\mathbb Z}}    

\newcommand{\calD}{\mathcal D}         
         
\newcommand{\calF}{\mathcal F}         
\newcommand{\calG}{\mathcal G}

\newcommand{\calK}{\mathcal K}         
\newcommand{\calL}{\mathcal L}         
         
\newcommand{\calN}{\mathcal N}         
         
\newcommand{\calP}{\mathcal P}

\newcommand{\calX}{\mathcal X}         
\newcommand{\calY}{\mathcal Y}         
\newcommand{\calZ}{\mathcal Z}         

\newcommand{\la}{\lambda}

\newcommand{\vphi}{\varphi}

\newcommand{\bmat}[4]{\begin{bmatrix}#1 & #2\\#3 & #4\end{bmatrix}}

\newcommand{\linspan}{\operatorname{span}}

\newcommand{\dd}{\text{d}}
\newcommand{\bbK}{\mathbf{K}}
\newcommand{\sgn}{\operatorname{sgn}}

\numberwithin{equation}{section}
\allowdisplaybreaks

\title[$L^\infty$-error bounds for approximations of the Koopman operator by kEDMD]{%
	$L^\infty$-error bounds for approximations of the Koopman operator by kernel extended dynamic mode decomposition
}
\author{Frederik Köhne$^1$}\address{$^1$Department of Mathematics, University of Bayreuth.\\ Mail: \textsc{\{frederik.koehne, anton.schiela\}@uni-bayreuth.de}}
\author{Friedrich M.\ Philipp$^2$}\address{$^2$Institute of Mathematics, Technische Universität Ilmenau, Germany.\\ Mail: \textsc{\{friedrich.philipp,manuel.schaller,karl.worthmann\}@tu-ilmenau.de}}
\author{Manuel Schaller$^2$}
\author{Anton Schiela$^1$}
\author{Karl Worthmann$^2$}

\thanks{F.P.\ was funded by the Carl Zeiss Foundation within the project DeepTurb--Deep Learning in and from Turbulence. He was further supported by the free state of Thuringia and the German Federal Ministry of Education and Research (BMBF) within the project THInKI--Th\"uringer Hochschulinitiative für KI im Studium. K.W.\ gratefully acknowledges funding by the Deutsche Forschungsgemeinschaft (DFG, German Research Foundation) -- Project-ID 507037103.}

\begin{document}
\begin{abstract}
Extended dynamic mode decomposition (EDMD) is a well-established method to generate a data-driven approximation of the Koopman operator for analysis and prediction of nonlinear dynamical systems. Recently, kernel EDMD (kEDMD) has gained popularity due to its ability to resolve the challenging task of choosing a suitable dictionary by using the kernel's canonical features and, thus, data-informed observables. In this paper, we provide the first pointwise bounds on the approximation error of kEDMD. The main idea consists of two steps. First, we show that the reproducing kernel Hilbert spaces of Wendland functions are invariant under the Koopman operator. Second, exploiting that the learning problem given by regression in the native norm can be recast as an interpolation problem, we prove our novel error bounds by using interpolation estimates. Finally, we validate our findings with numerical experiments.
\end{abstract}

\maketitle


\smallskip
\noindent \textbf{Mathematics subject classications (2020).} 37M99, 
41A05, 
47B32, 
47B33, 
65D12

\bigskip
\section{Introduction}
Introduced by B.O.\ Koopman in the 1930's \cite{Koop31}, the Koopman operator offers a powerful theoretical framework for data-driven analysis, prediction, and control of dynamical systems. Since the seminal paper~\cite{Mezi05} by I.\ Mezi{\'c} and the rapid development and success of deep learning techniques, Koopman's idea and related approaches to the learning of highly-nonlinear dynamics have received significant interest, see, e.g., \cite{BrunBudi22} and the references therein. In essence, the Koopman operator~$\calK_A$ predicts the dynamical behavior of the system through the lens of observable functions~$f : Y\to\R$, i.e., $\mathcal{K}_A f = f \circ A$, where the continuous map $A : X \to Y$ represents the flow map of a dynamical system on topological spaces $X$ and $Y$. Correspondingly, the nonlinear dynamics~$A$ are lifted into an infinite-dimensional function space, on which the Koopman operator~$\mathcal{K}_A$ acts linearly. Note that $\calK_A$ maps functions on $Y$ to functions on $X$.

Typically, data-driven techniques aim at learning a compression, i.e., a finite section of the Koopman operator, by evaluations of observables from a pre-defined dictionary of functions at a set of data points in state space. 
The most prominent representative is Extended Dynamic Mode Decompostion (EDMD), which has been successfully applied to climate prediction~\cite{AzenLin20}, molecular dynamics~\cite{WuNusk17}, turbulent flows~\cite{GianKolc18,Mezi13}, neuroscience~\cite{BrunJohn16}, and deep learning~\cite{DogrRedm20}, to name just a few.
However, as the qualitative and quantitative insight obtained by EDMD strongly depends on the dictionary, its choice is a central and delicate task.
An appealing way to resolve 
this problem is kernel EDMD (kEDMD; \cite{WillRowl15,klus2020kernel}), where the dictionary consists of the canonical features of an a-priori chosen kernel function centered at the data points. 

Despite the success of these Koopman-based methods and the plurality of works in the field, there are only a few results available concerning error bounds on the corresponding data-driven approximations in terms of the data. 
For deterministic systems, the first error bounds on EDMD were derived in~\cite{Mezi22} for ergodic sampling, followed by results in~\cite{ZhanZuaz23} for i.i.d.\ sampling. 
The first error bounds for control and stochastic systems were proven in~\cite{NuskPeit23} for both sampling regimes. 
Under significantly weaker conditions, novel error bounds --~ which also apply to discrete-time systems~-- were only recently shown in~\cite{PhilScha24}. 
Bounds for ResDMD, a variant of EDMD tailored to the extraction of spectral information, can be found in~\cite{ColbLi24}. 
Concerning kEDMD, the papers \cite{PhilScha23b} and \cite{PhilScha23} provide error bounds for both prediction and control using either i.i.d.\ or ergodic sampling.
In all of the mentioned results, the approximation error is split up into its two sources. 
The first one constitutes an estimation error resulting from finitely many data points. 
If the data is drawn randomly, the respective error is, as to be expected, of probabilistic nature. 
The second source quantifies the projection error onto the dictionary, and as such, heavily depends on the choice of the dictionary. 
When choosing the dictionary as a subset of an orthonormal basis, one may show that this projection error vanishes as the dictionary size tends to infinity~\cite{KordMezi18}. 
However, in applications, quantitative error estimates for finite dictionaries are key. 
In this context, for \textit{classical} EDMD schemes, only $L^2$-type bounds are available~\cite{ZhanZuaz23,SchaWort23} based on a dictionary of finite-element functions. 

In the recent preprint \cite{yadav2024approximation}, an alternative to EDMD is proposed for approximating the Koopman operator. This approach is based on approximation by uni- 
and multivariate Bernstein polynomials and inherits their property of uniform convergence. However, 
also the inherent restriction to the (perturbed) lattice structure of the Bernstein polynomials is currently preserved, see \cite[Subsection~5.1]{yadav2024approximation}.

In this paper, we provide the first uniform-error analysis for the well-established kernel EDMD method, which is very popular due to its flexiblility and efficiency. More precisely, we derive fully-deterministic bounds on the pointwise error of the kEDMD approximant~$\wh\calK_A$ of the Koopman operator $\calK_A$ on a native space (also called reproducing kernel Hilbert space; RKHS), $\calN(Y)$, generated by a kernel on $Y$. To this end, we first 
recall that regression in the native norm corresponds to an interpolation problem 
and show that $\wh\calK_A = S_\calX\calK_A$ holds, where $S_\calX$ denotes an orthogonal projection, based on the interpolation on a set of nodes $\calX\subset X$.
Second, we prove a bound on the full uniform 
approximation error $\|\calK_A f - \wh\calK_A f\|_\infty$ for $f\in\calN(Y)$, cf.\ \cref{thm:err_est}.  In this context, we require two main ingredients:
\begin{enumerate}
    \item [(i)] 
        The {\em inclusion property} $\calK_A\calN(Y)\subset\calN(X)$. If this holds, the Koopman operator $\calK_A : \calN(Y)\to\calN(X)$ is bounded, see \cref{lem:koopman-is-continuous-on-native-spaces}.
        Despite the negative result from~\cite{GonzAbud23}, which shows that in the case of {\em Gaussian} kernels the inclusion property can only hold for affine-linear dynamics~$A$, we prove in \cref{thm:koopman-well-defined} that the Koopman operator preserves Sobolev regularity, and hence $\calK_A\calN(Y)\subset\calN(X)$ holds for the Wendland native spaces \cite{Wend04}.
    \item [(ii)] 
        The uniform interpolation error in $\calN(X)$, i.e., the quantity $\|I - S_\calX\|_{\calN(X)\to C_b(X)}$, which can be controlled via sophisticated results from approximation theory~\cite{Fass05,MaddScha21,Wend04}.
\end{enumerate}
However, the calculation of $\wh\calK_A f$ for an observable $f$ requires the evaluation of $f$ at the $A(x_i)$, which might be prohibited by certain restrictions, e.g., when the dynamics are inaccessible at runtime.
We therefore analyze a variant of the kEDMD approximant defined by $\wh\calK_A^\calY = S_\calX\calK_A S_\calY$ and provide, again, bounds on the error in the uniform norm. 
Here, $\calY\subset Y$ is a second set of nodes, which may be chosen independently of $\calX$. This provides additional flexibility w.r.t.\ the kEDMD implementation and avoids the issue mentioned above.

The paper is organized as follows. In \cref{sec:native}, we recall the basic properties of RKHSs as well as interpolation and regression in these. 
Next, we build our abstract framework around the Koopman operator and kEDMD in \cref{sec:koopman} recalling that Koopman regression in native spaces by means of kEDMD is interpolation. The discussion is complemented by a characterization of the inclusion property in \cref{s:incprop}. In \cref{thm:err_est}, we prove uniform error bounds, which depend on the quantities $\|\calK_A\|_{\calN(Y)\to\calN(X)}$ and $\|I - S_\calX\|_{\calN(X)\to C_b(X)}$ mentioned above. In \cref{sec:wendland}, we introduce the Wendland radial basis functions (RBFs) and prove in \cref{thm:koopman-well-defined} that the Koopman operator preserves Sobolev regularity, which is then used in \cref{sec:uniform-error-bounds} to prove the uniform error bounds, cf.\ \cref{thm:main}. Finally, we provide numerical examples in \cref{sec:numerics} before conclusions are drawn in \cref{sec:conclusions}.

\subsection*{Notation}
We let $\N := \{1,2,3,\ldots\}$ and $\N_0 := \N\cup\{0\}$. For $m,n\in\N$, $m\le n$, we use the notation $[m:n] := \{k\in\N_0 : m\le k\le n\}$.

\section{Native spaces of kernels}\label{sec:native}
Let $X$ be a topological space, and let $k : X\times X\to\R$ be a bounded and continuous symmetric positive definite kernel function, $k\neq 0$. We adhere to the existing literature on kernel-based methods, where {\em positive definite} means that for all $n\in\N$, all $\zeta_1,\ldots,\zeta_n\in\R$, and all $x_1,\ldots,x_n\in X$ we have
\[
\sum_{i,j=1}^nk(x_i,x_j)\zeta_i\zeta_j\,\ge\,0.
\]
In other words, the symmetric kernel matrix $(k(x_i,x_j))_{i,j=1}^n$ is positive \emph{semi}-definite. For $z\in X$, define the function $\Phi_z : X\to\R$ by
\[
\Phi_x(z) := k(x,z),\qquad x\in X.
\]
The $\Phi_x$ are called the {\em canonical features} of $k$. We call $k$ {\em strictly} positive definite if the canonical features $\Phi_{x_1},\ldots,\Phi_{x_n}$ corresponding to any finite number of pairwise distinct points $x_1,\ldots,x_n\in\R^d$ are linearly independent. This is easily seen to be equivalent to the positive definiteness of the kernel matrices $(k(x_i,x_j))_{i,j=1}^n$. We refer to \cite[Theorem 3.6]{PaulRagh16} for further equivalent conditions.

It is well known \cite{PaulRagh16,SteiChri08} that the kernel function $k$ induces a unique Hilbert space $(\calN(X),\<\cdot\,,\cdot\>_{\calN(X)})$ of functions on $X$ such that $\Phi_x\in\calN(X)$ for all $x\in X$ and
\begin{align}\label{e:repr_pr}
f(x) = \<f,\Phi_x\>_{\calN(X)},\qquad x\in X ,\, f \in \calN(X).
\end{align}
In particular, the function evaluation $\delta_x : \calN(X)\to\R$, $\delta_x(f) = f(x)$, is a continuous linear functional on $\calN(X)$ for every $x\in X$. The Hilbert space $\calN(X)$ is called the {\em reproducing kernel Hilbert space} or the {\em native space} associated with $k$. In converse, if $H$ is any Hilbert space of functions on $X$, such that $\delta_x\in H^*$ for all $x\in X$, then $k(x,y):=\langle \delta_x,\delta_y\rangle_{H^*}$ defines a kernel function, such that $H$ is the corresponding native space.

From \cref{e:repr_pr} it is easily seen that the functions in $\calN(X)$ are also continuous and bounded, i.e., $\calN(X)\subset C_b(X)$. In fact, the respective embedding is continuous: for $f\in \calN(X)$ we have
\begin{align}\label{eq:norm_vs_norm}
\|f\|_{\calN(\domX)}=\sup_{g\in \calN(\domX)} \frac{\langle g,f\rangle_{\calN(\domX)}}{\|g\|_{\calN(\domX)}}\ge \sup_{x\in \domX} \frac{|f(x)|}{\|\Phi_x\|_{\calN(\domX)}} = \sup_{x\in\domX}\frac{|f(x)|}{\sqrt{k(x,x)}} \ge \|k\|_\infty^{-1/2}\|f\|_\infty.
\end{align}
In other words: the native norm is stronger than the $\sup$-norm.

The majority of the content in the following two subsections is well known and can be found in, e.g., \cite[Chapter 3]{PaulRagh16}. We have included it in order to be self-contained on the one hand and to introduce our notation on the other.

\subsection{Interpolation in native spaces}
Let $k$ be strictly positive definite. For a set of pair\-wise distinct \emph{centers} $\calX = \{x_i\}_{i=1}^N\subset X$, define the subspace
\begin{align*}
V_\calX = \spn\{\Phi_{x_i}: i\in[1:N]\}\subset\calN(X).
\end{align*}
By construction, the set $\{ \Phi_{x_i} : i\in [1:N]\}$ is a basis of 
$V_\calX$, which we call the 
\emph{canonical basis}. 
Given $\calZ = \{z_j\}_{j=1}^M\subset X$, define the (in general, rectangular) matrix $\bbK_{\calZ,\calX}\in\mathbb R^{M\times N}$ by
\[
(\bbK_{\calZ,\calX})_{ij}\coloneqq k(z_i,x_j) =  \Phi_{x_j}(z_i), \qquad i\in [1:M],\,j\in [1:N].
\]
If $\alpha_f\in\R^N$ is a basis representation of $f\in V_\calX$, i.e., $f = \sum_{j=1}^N(\alpha_{f})_j\Phi_{x_j}$, then $\bbK_{\calZ,\calX}\,\alpha_f$ consists of the values of $f$ at the points $z_i$, $i\in [1:N]$:
\[
\big(\bbK_{\calZ,\calX}\,\alpha_f\big)_i = \sum_{j=1}^N(\bbK_{\calZ,\calX})_{ij}(\alpha_{f})_j = \bigg(\sum_{j=1}^N(\alpha_{f})_j\Phi_{x_j}\bigg)(z_i) = f(z_i).
\]
Hence,
\begin{equation}\label{eq:transform}
\bbK_{\calZ,\calX}\,\alpha_f = f_\calZ,
\end{equation}
where $f_\calZ := (f(z_1),\ldots,f(z_M))^\top$. If we choose $M=N$ and $\calZ = \calX$, then, by strict positive definiteness of the kernel $k$, the kernel matrix $\bbK_{\calX} := \bbK_{\calX,\calX}\in\R^{N\times N}$ is symmetric positive definite and thus invertible. In this case, the interpolation problem:
\[
\text{Given $y\in\R^N$, find } f\in V_\calX\text{ such that }f(x_i) = y_i\text{ for all $i\in [1:N]$}
\]
possesses a unique solution $f\in V_\calX$, whose coefficients are given by
\[
\alpha_f = \bbK_{\calX}^{-1}f_\calX.
\]
We call $f$ the {\em $\calX$-interpolant} of $y$. We also call $f\in V_\calX$ the $\calX$-interpolant of any $g\in\calN(X)$ for which $g_\calX = f_\calX$. The vector $f_\calX = y$ consists of the basis coefficients of $f\in V_\calX$ in the {\em Lagrange basis}\footnote{In \cite{PaulRagh16}, the Lagrange basis is called the canonical partition of unity for $\calX$.} $\{\Phi_{x_i}^L\}_{i=1}^N$ of $V_\calX$, where $\Phi_{x_i}^L(x_j) = \delta_{ij}$ for $i,j\in [1:N]$: 
\[
  f = \sum_{i=1}^N (f_\calX)_i \Phi^L_{x_i} \quad \text{ for } \quad f \in V_\calX. 
\]
Thus, $\bbK_{\calX}$ also represents a change of basis from the canonical basis to the Lagrange basis, and solving the interpolation problem can be interpreted as sampling $f$ at the centers, followed by a change of basis. Finally, the native inner product on $V_\calX$ can be represented via basis representations, using \cref{e:repr_pr}:
\begin{equation}\label{eq:native_basis}
\langle f,g\rangle_{\calN(X)} = \bigg\<f,\sum_{j=1}^N\alpha_g \Phi_{x_j}\bigg\>_{\calN(X)} = \sum_{j=1}^N \alpha_g f(x_j) = \alpha_g^\top f_\calX,\qquad f,g\in V_\calX.
\end{equation}

\subsection{Regression in native spaces}\label{subsec:regression}
Let $f \in \mathcal N(\domX)$. It is our goal to find a good approximation $g\in V_\calX$ to $f$. If ``good approximation'' means ``best approximation in the native norm'', the corresponding regression problem reads
\begin{align}\label{eq:regression}
\min_{g\in V_\calX} \|f-g\|_{\mathcal N(\domX)}^2.
\end{align}
Obviously, the solution of \cref{eq:regression} is given by $g = S_\calX f$, where $S_\calX : \mathcal N(\domX)\to V_\calX$ denotes the orthogonal projection in $\mathcal N(\domX)$ onto $V_\calX$ with respect to the native inner product. On the other hand, note that \cref{eq:regression} is equivalent to
\[
0 = \langle f-g,\Phi_{x_j}\rangle_{\mathcal N(\domX)} = f(x_j) - g(x_j), \qquad j\in [1:N].
\]
Hence, {\it the solution $g\in V_\calX$ of \cref{eq:regression} coincides with the $\calX$-interpolant of $f_\calX$}.
As seen above, in the Lagrange basis the solution $g\in V_\calX$ has the basis coefficients $\alpha_g^L = g_\calX = f_\calX$, and thus
\begin{align}\label{eq:SxLagrange}
S_\calX f = \sum_{j=1}^N\alpha_{gj}^L\Phi^L_{x_j}=\sum_{j=1}^N f(x_j)\Phi^L_{x_j},\qquad f\in\calN(X).
\end{align}
In the canonical basis we have $\alpha_g = \bbK_{\calX}^{-1}g_\calX = \bbK_{\calX}^{-1}f_\calX$. Hence, the projection $S_\calX$ has the following representation in the canonical basis: 
\begin{align}\label{eq:Sx}
S_\calX f = \sum_{j=1}^N\big(\bbK_{\calX}^{-1}f_\calX\big)_j\Phi_{x_j},\qquad f\in\calN(X).
\end{align}

\begin{proposition}
The projection operator $S_\calX : \calN(X)\to V_\calX$ admits a bounded extension $\wt S_\calX : C_b(X)\to V_\calX$, which is given by
\[
\wt S_\calX f = \sum_{j=1}^N\big(\bbK_{\calX}^{-1}f_\calX\big)_j\Phi_{x_j}.
\]
\end{proposition}
\begin{proof}
The defined operator is certainly a well defined extension of $S_\calX$ to $C_b(X)$. For $f\in C_b(X)$ define $f_\calX = (f(x_j))_{j=1}^N=(\wt S_\calX f(x_j))_{j=1}^N$. To compute the native norm, we use \cref{eq:native_basis} to observe
\begin{align*}
\|\wt S_\calX f\|^2_{\calN(X)} = \alpha_{\wt S_\calX f}^\top f_\calX\le \|\alpha_{\wt S_\calX f}\|_{\ell^1}\|f_\calX\|_{\ell^\infty} = \|\bbK_{\calX}^{-1}f_\calX\|_{\ell^1}\|f_\calX\|_{\ell^\infty}\le \|\bbK_{\calX}^{-1}\|_{\ell^\infty\to\ell^1}\|f_\calX\|^2_{\ell^\infty},
\end{align*}
implying $\|\wt S_\calX  f\|_{\calN(X)}\le \|\bbK_{\calX}^{-1}\|^{1/2}_{\ell^\infty\to\ell^1}\|f\|_\infty$.
\end{proof}

\begin{remark}
As can be seen from the proof, the operator norm of $\wt S_\calX$ depends on the number and locations of the centers $x_i$ and we cannot expect a uniform bound with respect to these choices. 
\end{remark}

\section{Koopman approximation by regression in the native norm}\label{sec:koopman}
In this section, we develop our 
framework to construct approximations of 
the Koopman operator. Therein, the intimate relation between interpolation and regression in the native norm is a key 
ingredient to prove pointwise error bounds, as provided in the main result, \cref{thm:main} below.

\subsection{The Koopman operator on native spaces}
Let $X$ and $Y$ be topological spaces, and let $A : X\to Y$ be a continuous map. In the following, $A$ will model the flow map of a dynamical system for a fixed time step. Then $A$ induces the linear \emph{Koopman operator} $\calK_A : C_b(Y)\to C_b(X)$, which is defined by
\begin{align*}
\calK_A f := f\circ A,\qquad f\in C_b(Y).
\end{align*}
Note that the Koopman operator $\calK_A$ is a contraction, which follows directly from its definition. In fact, by noting $\calK_A\one_Y = \one_X$, we immediately see that $\|\calK_A\|_{C_b(Y)\to C_b(X)} = 1$.

Although the Koopman operator $\calK_A$ is a linear object which captures the full information on the dynamics $A$, in real-world applications it is usually unknown, and Koopman-based data-driven methods aim at learning the operator $\calK_A$ from finitely many observations of the dynamics $A$ at a finite number of states $x_i\in X$, i.e., from $\calK_A f_j(x_i) = f_j(A(x_i))$, where the $f_j\in C_b(Y)$ are called {\em observables}.

In what follows, we consider particular observables in $\calN(Y)$. To this end, let $k_X$ and $k_Y$ be continuous and bounded strictly positive-definite kernel functions on $X$ and $Y$ with canonical features $\Phi_x$ and $\Psi_y$ for $x\in X$ and $y\in Y$, respectively.
We shall assume occasionally that $\calK_A$ maps $\calN(Y)$ into $\calN(X)$, i.e.,
\begin{align}\label{eq:invariant}
\calK_A\,\calN(Y)\subset\calN(X).    
\end{align}
We  will see in \cref{sec:wendland} that this condition is satisfied for appropriately chosen Wendland functions under mild conditions on the domains and on $A$.

If \cref{eq:invariant} holds, we observe the following adjoint property of $\koopman$ on native spaces
\begin{equation}\label{eq:koopmanadjoint}
\langle \koopman f,\Phi_x\rangle_{\calN(X)}=f(A(x))=\langle f,\Psi_{A(x)}\rangle_{\calN(Y)}, \qquad f\in\calN(Y),\,x\in X.
\end{equation}

\begin{lemma}\label{lem:koopman-is-continuous-on-native-spaces}
If \eqref{eq:invariant} holds, then the restriction $\calK_A|_{\calN(Y)} : \calN(Y)\to\calN(X)$ is a bounded operator.
\end{lemma}
\begin{proof}
Suppose that \eqref{eq:invariant} holds. By the closed graph theorem, for the boundedness of $\calK_A : \calN(Y)\to\calN(X)$, it suffices to show that this operator is closed. To this end, let $f,f_n\in\calN(Y)$ and $g\in\calN(X)$ such that $\|f_n-f\|_{\calN(Y)}\to 0$ and $\|\calK_A f_n - g\|_{\calN(X)}\to 0$ as $n\to\infty$. We have to show that $g = \calK_A f$. To this end, for any $x\in X$, we compute via \cref{eq:koopmanadjoint}
\begin{align*}
g(x)
&= \<g,\Phi_x\>_{\calN(X)} = \lim_{n\to\infty}\<\koopman f_n,\Phi_x\>_{\calN(X)}  = \lim_{n\to\infty}\<f_n,\Psi_{A(x)}\>_{\calN(Y)}\\
&= \<f,\Psi_{A(x)}\>_{\calN(Y)} = f(A(x)) = (\calK_Af)(x),
\end{align*}
and hence obtain $g = \calK_A f$. 
\end{proof}

Note that \cref{lem:koopman-is-continuous-on-native-spaces} still holds 
if we drop the continuity and the strict positive definiteness of the kernels $k_X$ and $k_Y$. In \cref{s:incprop}, we characterize the inclusion property \eqref{eq:invariant} and infer that typically $\|\calK_A\|_{\calN(Y)\to\calN(X)} > 1$.

\subsection{Kernel EDMD}\label{subsec:kEDMD}
Extended Dynamic Mode Decomposition (EDMD) is a data-driven method which aims at approximating the Koopman operator~$\calK_A$ based upon evaluations of finite sets of observable functions $\calG\subset C_b(X)$, $\calF\subset C_b(Y)$, frequently called {\em dictionaries}\footnote{
Often in the literature, only the case $X = Y$ and $\calG = \calF$ is discussed.}, at nodes $x_1,\ldots,x_m\in X$ and the corresponding 
values $A(x_1),\ldots A(x_m)$, respectively. Setting $V := \linspan\calG$ and $W := \linspan\calF$, EDMD builds a data matrix $\bbK$ representing an operator $\calK_A' : W\to V$ which approximates the compression $P_V\calK_A|_W$ of the Koopman operator, regarded as a map between some $L^2$-spaces over $X$ and $Y$, respectively. Here, $P_V$ and $P_W$ denote the $L^2$-orthogonal projections onto $V$ and $W$, respectively. Having approximated the finite-dimensional compression, one then has to consider the projection error corresponding to the difference of $P_V \calK_A|_W$ and $\calK_A|_W$ or, if one wants to also predict observables not contained in the dictionary $\calF$,  to the difference between $P_V \calK_A P_W$ and $\calK_A$. Clearly, this projection error strongly depends on the chosen dictionaries $\calF$ and $\calG$. Hence, the choice of a suitable dictionary is a central and delicate task.

In kernel EDMD (kEDMD; \cite{WillRowl15}), this issue
is alleviated by means of data-driven dictionaries. Therein, one has $N=m$, i.e., the dictionary size coincides with the amount of data points, and one adds a set of centers $\calY = \{y_j\}_{j=1}^M$ in $Y$ to define the corresponding dictionaries as $\calF = \{\Psi_{y_i} : i\in [1:M]\}$ and $\calG = \{\Phi_{x_j} : j\in [1:N]\}$. Hence, the problem of choosing the dictionaries reduces to the mere choice of a kernel. Note that, in this case, $V = V_\calX$ and $W = V_\calY$ according to \cref{sec:native}. Moreover, $M$ and the centers $y_j$ can be freely chosen. For example, we may choose $M=N$ and $y_i = A(x_i)$ or $y_i = x_i$ if $A : X\to X$ (i.e., $Y=X$).

In this setting, the EDMD matrix approximant $\bbK$ reads
\begin{align}\label{eq:kedmd_mat}
\bbK = \bbK_{\calX}^{-1}\bbK_{A(\calX),\calY}.
\end{align}

\begin{proposition}\label{prop:mat_repr}
The matrix $\bbK$ represents the operator $\wt S_\calX\calK_A|_{V_\calY}$ with respect to the canonical bases of $V_\calX$ and $V_\calY$.
\end{proposition}
\begin{proof}
For $j\in [1:M]$ we have
\begin{align*}
\wt S_\calX\calK_A\Psi_{y_j} = \wt S_\calX(\Psi_{y_j}\circ A) = \sum_{i=1}^N\big(\bbK_{\calX}^{-1}[\Psi_{y_j}\circ A]_\calX\big)_i\Phi_{x_i}.
\end{align*}
Hence, the $j$-th column of the matrix representing $\wt S_\calX\calK_A|_{V_\calY}$ is given by
\[
\bbK_{\calX}^{-1}[\Psi_{y_j}\circ A]_\calX = \bbK_{\calX}^{-1}\big(k_Y(A(x_i),y_j)\big)_{i=1}^N = \bbK_{\calX}^{-1}\bbK_{A(\calX),\calY}e_i,
\]
where $e_i$ denotes the $i$-th standard basis vector. The claim now follows from \cref{eq:kedmd_mat}.
\end{proof}

Fixing $\calX$, a natural approximant of the full operator $\calK_A|_{\calN(Y)}$ is therefore given by
\begin{align}\label{eq:hatKY}
\wh\calK_A^{\calY} := \wt S_\calX\calK_A S_\calY.
\end{align}
For the special choice $\calY = A(\calX)$ we set
\begin{align}\label{eq:hatK}
\wh\calK_A := \wh\calK_A^{A(\calX)} = \wt S_\calX\calK_A S_{A(\calX)}.
\end{align}

\begin{proposition}\label{prop:methods_equal}
We have $\wh\calK_A = \wt S_\calX\calK_A$. If~\cref{eq:invariant} holds, then $\wh\calK_A^\calY = S_\calX\calK_A$ if and only if $A(\calX)\subset\calY$.
\end{proposition}
\begin{proof}
We have to show that $\wt S_\calX\calK_A S_{A(\calX)} = \wt S_\calX\calK_A$, i.e., $\wt S_\calX\calK_A(I - S_{A(\calX)}) = 0$. If $f\in\calN(Y)$, then $g = (I - S_{A(\calX)})f$ satisfies $g(A(x_i)) = 0$ for $i\in[1:N]$. And since $\wt S_\calX\calK_A g$ is the $\calX$-interpolant of $g\circ A$, it is uniquely determined by $(g\circ A)(x_i) = 0$ and thus vanishes.

Assume that \cref{eq:invariant} holds. By strict positive definiteness of $k_Y$, $A(\calX)\subset \calY$ is equivalent to $V_{A(\calX)}\subset V_\calY$, which is in turn equivalent to $\langle S_\calY f,v\rangle_{\calN(Y)} = \langle f,v\rangle_{\calN(Y)}$ for all $f\in\calN(Y)$ and all $v\in V_{A(\calX)}$, since $S_\calY$ is the orthogonal projection onto $V_\calY$. Finally, this is equivalent to $\wh\calK_A^\calY = S_\calX\calK_A$, as the following computation shows:
\begin{align*}
S_\calX \koopman S_\calY f =  S_\calX \koopman f \quad &\Longleftrightarrow \quad 
\langle \koopman S_\calY f,\Phi_{x_i}\rangle_{\calN(X)}=\langle \koopman f,\Phi_{x_i}\rangle_{\calN(X)}\\
&\Longleftrightarrow \quad 
\langle S_\calY f,\Psi_{A(x_i)}\rangle_{\calN(Y)}=\langle f,\Psi_{A(x_i)}\rangle_{\calN(Y)},
\end{align*}
where we have used \cref{eq:koopmanadjoint}.
\end{proof}

We are now able to provide deterministic error bounds in the uniform norm for the approximants $\wh\calK_A^\calY f$ and $\wh\calK_A f = \wh\calK_A^{A(\calX)}f$ of $\calK_A f$ for observables $f\in\calN(Y)$. These in particular imply pointwise bounds on the approximation error of data-driven predictions via kEDMD. 

\begin{theorem}\label{thm:err_est}
If the inclusion $\calK_A \calN(Y)\subset \calN(X)$ holds \braces{cf.\ \eqref{eq:invariant}}, then for the Koopman approximation $\wh \calK_A$ defined in \cref{eq:hatK} we have the error bound
\[
\big\|\calK_A - \wh\calK_A\big\|_{\calN(Y)\to C_b(X)}\,\le\,\|I - S_\calX\|_{\calN(X)\to C_b(X)}\|\calK_A\|_{\calN(Y)\to\calN(X)}.
\]
and, for the error of $\wh\calK_A^\calY$ defined in~\cref{eq:hatKY} it holds that
\begin{align*}
\|\calK_A - \wh\calK_A^\calY\|_{\calN(Y)\to C_b(X)}
\le \|I-S_\calY\|_{\calN(Y)\to C_b(Y)} + \|I-S_\calX\|_{\calN(X)\to C_b(X)}\|\calK_A\|_{V_\calY\to \calN(X)}.
\end{align*}
In particular, for $f\in\calN(Y)$ we have the uniform bounds on the approximation error
\[
\big\|\calK_Af - \wh\calK_Af\big\|_\infty\,\le\,\|I - S_\calX\|_{\calN(X)\to C_b(X)}\|\calK_A\|_{\calN(Y)\to\calN(X)}\|f\|_{\calN(Y)}.
\]
and
\begin{align*}
\|\calK_Af - \wh\calK_A^\calY f\|_\infty\le \|f-S_\calY f\|_\infty + \|I-S_\calX\|_{\calN(X)\to C_b(X)}\|\calK_A\|_{V_\calY\to\calN(X)}\|S_\calY f\|_{\calN(Y)}.
\end{align*}
\end{theorem}
\begin{proof}
We have $\calK_A - \wh\calK_A = \calK_A - S_\calX\calK_A = (I - S_\calX)\calK_A$, from which the first estimate readily follows. On the other hand, concerning the error caused by the approximant $\wh\calK_A^\calY$, we observe that $\calK_A - \wh\calK_A^\calY = \calK_A - S_\calX\calK_A S_\calY = \calK_A (I-S_\calY) + (I-S_\calX)\calK_A S_\calY$, so that
\begin{align*}
\big\|\calK_A - \wh\calK_A^\calY\big\|_{\calN(Y)\to C_b(X)}
&\le \|\calK_A (I-S_\calY)\|_{\calN(Y)\to C_b(X)} + \|(I-S_\calX)\calK_A S_\calY
\|_{\calN(Y)\to C_b(X)}\\
&\le \|\calK_A\|_{C_b(Y)\to C_b(X)}\|I-S_\calY\|_{\calN(Y)\to C_b(Y)}\\
&\hspace*{2cm}+ \|I-S_\calX\|_{\calN(X)\to C_b(X)}\|\calK_A\|_{V_y\to\calN(X)}\|S_\calY\|_{\calN(Y)\to V_y}.
\end{align*}
The claim is now a consequence of $\|\calK_A\|_{C_b(Y)\to C_b(X)} = 1$ and $\|S_\calY\|_{\calN(Y)\to V_y} = 1$.
\end{proof}
The derived bounds feature two central components. The first one is the norm of the Koopman operator which will be bounded in \cref{sec:wendland} for native spaces of Wendland functions. The second ingredient are interpolation estimates, which, for various choices of the kernel, may be readily taken from the 
literature. They usually depend on a density measure on the $x_i$ in $X$, commonly called the \emph{fill distance}, cf.~\cref{thm:intest} of \cref{sec:uniform-error-bounds}. 

In the following remark, we 
provide further, more-detailed comments on the uniform error bounds presented in \cref{thm:err_est}.

\begin{remark}

\smallskip
(a) In \cite{PhilScha23b} (see also \cite{PhilScha23}), probabilistic bounds on the approximation error $\|\calK_A - \wh\calK_A\|_{\calN(Y)\to L^2(X,\mu)}$ were provided for the case where \cref{eq:invariant} holds. Here, the measure $\mu$ is a probability measure on $\calX$ and the $x_i$ are sampled randomly with respect to $\mu$. Of course, since $\|\cdot\|_{L^2(X,\mu)}\le\|\cdot\|_{C_b(X)}$, \cref{thm:err_est} also yields bounds on the approximation error in the $\|\cdot\|_{\calN(Y)\to L^2(X,\mu)}$ norm. However, the corresponding upper bounds are of a completely different nature compared to those in \cite{PhilScha23b,PhilScha23}. 

\smallskip
(b) For the estimate on $\wh\calK_A^\calY$, we do not actually need the inclusion \cref{eq:invariant}. Instead, one only has to require that $\calK_A\Psi_{y_j}\in\calN(X)$ for all $j\in[1:M]$. Then, as a linear operator on a finite-dimensional space, $\calK_A|_{V_\calY} : V_\calY\to\calN(X)$ is bounded. In the absence of \cref{eq:invariant}, this bound may, however, not be uniform in the number and location of the centers $y_i$.

\smallskip
(c) When applying $\wh\calK_A^\calY$, one might be tempted to use $\calY = \calX$, i.e., $\wh K_A^\calX = S_\calX\calK_A S_\calX$. However, this approach might cause comparably large 
errors at the boundary of $X$ if $X$ is not invariant under the flow $A$, i.e., $x_i\in X$, but $A(x_i)\in Y\backslash X$. An example for this behavior is provided in \cref{sec:numerics}.

\smallskip
(d) Similarly as in the proof of \cref{prop:inv_charac}, one can prove that
\begin{align}\label{e:sup2}
\|\calK_A\|_{V_\calY\to\calN(X)} = \sup_\calX\big\|\bbK_{\calX}^{-1/2}\bbK_{A(\calX),\calY}\bbK_{\calY}^{-1}\bbK_{\calY,A(\calX)}\bbK_{\calX}^{-1/2}\big\|_{2\to 2}^{1/2}
\end{align}
In particular, if $\calY = A(\calX)$, then $\|\calK_A\|_{V_\calY\to\calN(X)} = \|\calK_A\|_{\calN(Y)\to\calN(X)}$.
\end{remark}

\subsubsection*{Implementation aspects}
Although the two approximating operators $\wh\calK_A$ and $\wh\calK_A^\calY$ are mapping between the infinite-dimensional spaces $\calN(Y)$ and $\calN(X)$, their action can be computed with the help of certain kernel matrices. In what follows, we shall investigate to which extent the dynamics, given by $A : X\to Y$, has to be observed for the calculation of the approximants. We distinguish between observations $\Psi_{y_j}(A(x_i)) = k_Y(A(x_i),y_j)$ with the fixed kernel observables $\Psi_{y_j}$ and measurements $f(A(x_i))$ with arbitrary observables $f\in\calN(Y)$.

First of all, we note that the action of the operator $\wh\calK_A$ on an observable $f\in\calN(Y)$ requires the knowledge of the values $f(A(x_i))$. In fact, we have
\[
\wh\calK_Af = S_\calX\calK_A f = \sum_{j=1}^N\big(\bbK_{\calX}^{-1}[f\circ A]_\calX\big)_j\Phi_{x_j}, \quad\text{i.e., $\alpha_{\wh\calK_Af} = \bbK_{\calX}^{-1}[f\circ A]_\calX$}.
\]
If the values $f(A(x_i))$ are not available, one may resort to the variant $\wh\calK_A^\calY$ with a ``user-defined'' set of centers $\calY$ in $Y$. Indeed, we have
\[
\wh\calK_A^\calY f = S_\calX\calK_A S_\calY f = \sum_{j=1}^N\big(\bbK_{\calX}^{-1}[(S_\calY f)\circ A]_\calX\big)_j\Phi_{x_j},
\]
and since $S_\calY f(A(x_i)) = (\bbK_{A(\calX),\calY}\bbK_{\calY}^{-1}f_\calY)_i$, we obtain
\begin{equation}\label{eq:ktilde}
\alpha_{\wh\calK_A^\calY f} = \bbK_{\calX}^{-1}\bbK_{A(\calX),\calY}\bbK_{\calY}^{-1}f_\calY,
\end{equation}
which only requires the knowledge of the values $\Psi_{y_j}(A(x_i)) = k_Y(A(x_i),y_j)$, $i\in [1:N]$, $j\in [1:M]$, in the matrix $\bbK_{A(\calX),\calY}$ and evaluations of $f$ at the user-defined nodes $y_j$ which are independent of the dynamics. Hence, $\wh\calK_A^\calY f$ may be computed for any $f\in\calN(Y)$ by using the previously observed ``offline measurements'' $k_Y(A(x_i),y_j)$ and by evaluations of $f$ at the $y_j$ at runtime.

For explicit matrix representations we consider $\wh\calK_A^\calY|_{V_\calY} : V_\calY\to V_\calX$. Depending on the choice of bases we have
\begin{align*}
  \wh \calK^\calY_A|_{V_\calY}\sim \bbK_{\calX}^{-1}\bbK_{A(\calX),\calY} \qquad \text{ (canonical basis), }\\
  \wh \calK^\calY_A|_{V_\calY} \sim  \bbK_{A(\calX),\calY}\,\bbK_{\calY}^{-1} \qquad \text{ (Lagrange basis). } 
\end{align*}
The first representation is from \cref{prop:mat_repr} and the second follows by taking into account the appropriate basis transformation matrices between the canonical basis and the Lagrange basis. 
 The choice $\calY=\calX$ allows to regard $\wh \calK^\calX_A|_{V_\calX}: V_\calX\to V_\calX$ as a linear automorphism, and thus a meaningful computation of powers and of eigenvalues.

\subsection{Multi-step predictions}
If $Y=X$, and hence $\calK_A : C_b(X)\to C_b(X)$, powers $\calK_A^n$ of the Koopman operator are well defined, where $\calK_A^n f=f\circ A\circ \dots \circ A$, the composition taken $n$ times. In this setting, we would like to approximate $\calK_A^n$ for multi-step predictions. 

Let us assume that $\bbK_{\calX}^{-1}$ and $\bbK_{A(\calX),\calX}$  are available. Then multi-step predictions, i.e., approximations of $\calK_A^n f$ by $\wh\calK_A^n f$, only require the knowledge of 
the $f(A(x_i))$. In fact, we have
\[
\alpha_{\wh \calK_A^nf} = (\bbK_{\calX}^{-1}\bbK_{A(\calX),\calX})^{n-1}\bbK_{\calX}^{-1}[f\circ A]_\calX.
\]
This can be easily proved by induction and the fact that for any $h\in V_\calX$ we can evaluate $h\circ A$ at the points $x_i\in\calX$ by $[h\circ A]_\calX = \bbK_{A(\calX),\calX}\,\alpha_h$.

Once $\alpha_{\wh \calK_A^nf}$ has been computed, it is possible to evaluate $\wh \calK_A^nf$ at an arbitrary finite set of points $\calZ \subset X$ by the formula
\[
  [\wh \calK_A^n f]_\calZ = \bbK_{\calZ,\calX}\,\alpha_{\wh \calK_A^nf}.
\]
For multi-step predictions with $\wh\calK_A^\calX = S_\calX\calK_A S_\calX$ (i.e., $\calY = \calX$), the computation of 
\[
(\wh\calK_A^\calX)^n =(S_\calX\calK_A S_\calX)^n= \wh\calK_A^{n}S_\calX
\]
only requires the knowledge of $f(x_i)$ and we obtain the formula
\[
\alpha_{(\wh \calK_A^\calX)^nf} = (\bbK_{\calX}^{-1}\bbK_{A(\calX),\calX})^{n}\bbK_{\calX}^{-1}f_\calX.
\]
Our error bounds for single step predictions can be extended to multi-step predictions by a straightforward induction argument.

\begin{theorem}\label{thm:multistep}
If $X=Y$ and \cref{eq:invariant} holds, i.e., $\calK_A\calN(X)\subset\calN(X)$, then
\[
\|\calK_A^n-\wh \calK_A^n\|_{\calN(X)\to C_b(X)} \le \|I-S_\calX\|_{\calN(X)\to C_b(X)}\sum_{i=1}^n\|\calK_A\|^i_{\calN(X)\to \calN(X)}
\]
and 
\[
\|\calK_A^n-(\wh \calK_A^\calX)^n\|_{\calN(X)\to C_b(X)} \le \|I-S_\calX\|_{\calN(X)\to C_b(X)}\sum_{i=0}^n\|\calK_A\|^i_{\calN(X)\to \calN(X)}
\]
\end{theorem}
\begin{proof}
The statement is true for $n=1$ by \cref{thm:err_est}. Assume that it holds for $n-1$. We observe
\begin{align}\label{eq:dasding}
\calK_A^n-\wh \calK_A^n =\calK_A^n-\wh \calK_A^{n-1}S_\calX\calK_A=\calK_A^{n-1}(I-S_\calX)\calK_A+(\calK_A^{n-1}-\wh\calK_A^{n-1})S_\calX\calK_A. 
\end{align}  
Thus, using $\|S_\calX\|_{\calN(X)\to \calN(X)}=1$, $\|\calK_A\|_{C_b(X)\to C_b(X)}\le 1$, and the induction hypothesis,
\begin{align*}
\|\calK_A^n-\wh \calK_A^n\|_{\calN(X)\to C_b(X)}
&\le \|\calK_A^{n-1}\|_{C_b(X)\to C_b(X)}\|I-S_\calX\|_{\calN(X)\to C_b(X)}\|\calK_A\|_{\calN(X)\to \calN(X)}\\
&\hspace*{.5cm}+\|\calK_A^{n-1}-\wh \calK_A^{n-1}\|_{\calN(X)\to C_b(X)}\|S_\calX\calK_A\|_{\calN(X)\to \calN(X)}\\
&\le \|I-S_\calX\|_{\calN(X)\to C_b(X)}\|\calK_A\|_{\calN(X)\to \calN(X)}\\
&\hspace*{.5cm}+\|I-S_\calX\|_{\calN(X)\to C_b(X)}\left(\sum_{i=1}^{n-1}\|\calK_A\|^i_{\calN(X)\to \calN(X)}\right)
\|\calK_A\|_{\calN(X)\to \calN(X)}.
\end{align*}
This proves our first statement. The identity
\[
\calK_A^n-(\wh \calK_A^\calX)^n =\calK_A^n-\wh \calK_A^nS_\calX= \calK_A^n(I-S_\calX)+( \calK_A^n-\wh \calK_A^n)S_\calX.
\]
shows that the second statement follows from the first.
\end{proof}

\begin{remark}\label{r:multi_step_might _be_bad}
If $x\mapsto k(x,x)$ is constant on $X$, \cref{c:norm_large} implies that the operator norm $\|\calK_A\|_{\calN(X)\to\calN(X)} > 1$ (except the case $k_X = k_Y(A(\cdot),A(\cdot))$). Hence, for the right-hand sides of the estimates in \cref{thm:multistep} to be small, the interpolation error $\|I - S_\calX\|_{\calN(X)\to C_b(X)}$ might have to compensate for a potentially large number.
\end{remark}

\section{Koopman operators on Native Spaces of Wendland functions}\label{sec:wendland}
In this section, we show that the central invariance assumption~\cref{eq:invariant} of the error estimate provided in \cref{thm:err_est} is satisfied for native spaces of Wendland functions. The key ingredients are an intimate relation between these spaces with Sobolev spaces (\cref{subsec:wendland}) and that the Koopman operator preserves Sobolev regularity if the underlying flow has sufficient regularity (\cref{subsec:preserves}).

Throughout this section, let $X = \Omega_X$ and $Y = \Omega_Y$ be bounded Lipschitz domains in $\R^d$, i.e., they have a Lipschitz boundary, cf.\ \cref{app:boundary}, and let $A : \Omega_X\to\Omega_Y$ be a $C^1$-diffeomorphism such that
\begin{align}\label{e:bb}
\inf_{x\in\Omega_X}|\det DA(x)| > 0.
\end{align}
In the next subsection, we recall the compactly supported Wendland RBFs which induce the appropriate native spaces we shall be working with.

\subsection{Wendland functions and Sobolev spaces}\label{subsec:wendland}
For the reader's convenience, we briefly summarize the construction of the compactly supported, positive definite RBFs from \cite{Wend04}.
For $\ell \in \N$, define the \emph{base} function $\phi_\ell(r) = (1 - r)^\ell_+ = \max(1 - r ,0)^\ell$, $r \ge 0$.
For $r \ge 0$ and measurable functions $\phi: [0, \infty) \to [0, \infty)$ define
\[
\mathcal I \phi (r) = \int_r^\infty t \phi(t) dt. 
\]
Now, for given dimension $d\ge 1$ and smoothness $k\in\N_0$, set 
\[
\phi_{d,k} = \mathcal I ^k \phi_{\floor{\frac d2 } + k + 1},
\]
which is easily verified to be of the form 
\[
\phi_{d,k}(r)
= 
\begin{cases}
p_{d,k}(r), & 0 \le r \le 1 \\
0, & r > 1
\end{cases}
\]
with a univariate polynomial $p_{d,k}$ of degree $\floor{\frac d2} + 3k + 1$ and $\phi_{d,k}\in C^{2k}([0,\infty))$, see \cite[Theorem 9.13]{Wend04}. A recursive scheme for computing the coefficients of the polynomial $p_{d,k}$ can be found in \cite[Theorem 9.12]{Wend04}. Correspondingly, one may define the compactly supported RBF $\Phi_{d,k,x} : \R^d\to\R$ by
\[
\Phi_{d,k,x} := \phi_{d,k}(\|x - \cdot\|_2),
\]
where $d\in\N$, $k\in\N_0$, $x\in\R^d$ and set $\Phi_{d,k} := \Phi_{d,k,0}\in C^k(\R^d)$. The corresponding native space on a set $\Omega\subset\R^d$ will be denoted by $\calN_{\Phi_{d,k}}(\Omega)$. Note that the associated kernel function $k(x,y) = \phi_{d,k}(\|x-y\|_2)$ is bounded and continuous such that $\calN_{\Phi_{d,k}}(\Omega)\subset C_b(\Omega)$. Moreover, $k$ is strictly positive definite\footnote{Note that {\it positive definite} in \cite{Wend04} means strictly positive definite as defined in the beginning of \cref{sec:native}.} on $\R^d$ (and thus on any subset $\Omega\subset\R^d$) by \cite[Theorem 9.13]{Wend04}.

Given a bounded Lipschitz domain $\Omega\subset\R^d$, we denote the $L^2$-Sobolev space of regularity order $\sigma\ge 0$ on $\Omega$ by $H^\sigma(\Omega)$. Note that, as $\Omega$ has a Lipschitz boundary, the two predominant definitions of $H^\sigma(\Omega)$ in the literature via the Fourier transform on the one hand and via weak derivatives on the other coincide, see \cite[Theorem 3.30]{McLe00}. Therefore, we may use the following norm on $H^\sigma(\Omega)$ for $\sigma = \lfloor\sigma\rfloor + r$, cf.\ \cite[Chapter 3]{McLe00}:
\[
\norm{f}_{H^\sigma(\Omega)}^2 := \sum_{|\alpha|\le\lfloor\sigma\rfloor} \norm{D^\alpha f}_{L^2(\Omega)}^2 + \chi_{(0,1)}(r)\cdot\sum_{|\alpha|=\lfloor\sigma\rfloor}\int_\Omega\int_\Omega\frac{\big|D^\alpha f(x) - D^\alpha f(y)\big|^2}{\|x-y\|_2^{d+2r}}\,\dd x\,\dd y.
\]
The following theorem summarizes the connections between the RKHS corresponding to the RBFs defined above and Sobolev spaces.

\begin{theorem}\label{thm:native-spaces-sobolev-spaces}
Let $k\in\N_0$ and $d\in\N$ \braces{where $d\ge 3$ if $k=0$}. Let $\sigma_{d,k} = \frac{d+1}{2} + k$. Then
\[
\native _{\Phi_{d,k}}(\R^d) = H^{\sigma_{d,k}}(\R^d)
\]
with equivalent norms. If the bounded domain $\Omega\subset\R^d$ has a Lipschitz boundary, then
\[
\native _{\Phi_{d,k}}(\Omega) = H^{\sigma_{d,k}}(\Omega)
\]
with equivalent norms.
\end{theorem}
\begin{proof}
The first identity is shown in \cite[Theorem 10.35]{Wend04}. For odd $d$ (i.e., $\sigma_{d,k}\in\N$), the second follows directly from \cite[Corollary 10.48]{Wend04} when applied to the decay proven in \cite[Theorem 10.35]{Wend04}. The restriction to integer-valued Sobolev orders in \cite[Corollary 10.48]{Wend04} is not necessary. This restriction is due to the proof of \cite[Corollary 10.48]{Wend04}, where \cite[Theorem 1.4.5]{BrenScot94} is invoked, which proves the existence of a bounded extension operator $E : H^\sigma(\Omega)\to H^\sigma(\R^d)$ for $\sigma\in\N$. However, \cite[Theorem A.4]{McLe00} shows that such a bounded extension operator also exists for fractional Sobolev orders $\sigma\in\R$.
\end{proof}

\subsection{The Koopman operator preserves Sobolev regularity}\label{subsec:preserves}
In \cref{sec:koopman}, we introduced the Koopman operator as a linear contraction from $C_b(\Omega_Y)$ to $C_b(\Omega_X)$. The law $f\mapsto f\circ A$, however, makes sense for {\em any} function $f : Y\to\R$. In this part of the paper, we prove that the Koopman operator maps $L^2(\Omega_Y)$ boundedly into $L^2(\Omega_X)$ and, in general, Sobolev spaces over $\Omega_Y$ boundedly into Sobolev spaces over $\Omega_X$ and thus preserves the corresponding regularity---as long as the map $A$ has the same regularity. 

For a domain $\Omega\subset\R^d$, by $C_b^k(\Omega,\R^d)$ we denote the space of all $C^k(\Omega,\R^d)$-maps with bounded derivatives up to order $k$.

\begin{theorem}\label{thm:koopman-well-defined}
Assume in addition that $A\in C_b^m(\Omega_X,\R^d)$ for some $m \in \N$, $m > d/2$. Then for all $\sigma \le m$ the linear Koopman operator 
\begin{align}\label{eq:Koop_pres_reg}
\koopman : H^\sigma(\Omega_\domY) \to H^\sigma(\Omega_\domX)
\end{align}
is well defined and bounded. 
\end{theorem}
\begin{proof}
The proof consists of five steps: first of all, we prove the well-definedness of \cref{eq:Koop_pres_reg} for $H^0 = L^2$. In Step 2, we prove it for $H^1$ and generalize this result to $H^\sigma$, $\sigma\in\N_0$, in Step 3. The boundedness and fractional Sobolev spaces are treated in Steps 4 and 5.

\smallskip
\noindent {\bf Step~1}: \cref{eq:Koop_pres_reg} with $\sigma=0$. By the change of variables formula, for $f\in L^2(\Omega_Y)$ we have
\begin{align*}
\int_{\Omega_\domX} \abs{\koopman f(x)}^2\,\dd x = \int_{\Omega_X}\abs{f(A(x))}^2\,\dd x 
= \int_{\Omega_Y}\abs{f(y)}^2\abs{\det D\inv A(y)}\,\dd y.
\end{align*}
This proves that $\calK_A : L^2(\Omega_Y)\to L^2(\Omega_X)$ is well defined and bounded with
\[
\|\calK_A\|_{L^2(\Omega_Y)\to L^2(\Omega_X)}\le\|\det DA^{-1}\|_{C_b(\Omega_Y)}^{1/2}.
\]
{\bf Step 2.} \cref{eq:Koop_pres_reg} with $\sigma=1$. Let $f\in H^1(\Omega_Y)$. Then $\calK_A f\in L^2(\Omega_X)$ by Step 1. To prove $\calK_A f \in H^1(\Omega_\domX)$ it thus remains to show that $\koopman f$ possesses weak partial derivatives in $L^2(\Omega_X)$ up to order $1$. To see this, recall that $C^\infty(\Omega_Y)\cap H^1(\Omega_Y)$ is dense in $H^1(\Omega_Y)$ \cite[Theorem 3.17]{Adams2003}. Thus, we may consider a sequence $f_n \subset C^\infty(\Omega_Y)\cap H^1(\Omega_Y)$ with $f_n \to f$ as $n\to\infty$ in $H^1(\Omega_Y)$. Then $\norm{\koopman f_n - \koopman f}_{L^2(\Omega_X)} \to 0$ by Step 1 and, since $A$ is continuously differentiable, $\koopman f_n \in C^1(\Omega_X)$. In particular, by the chain rule we get
\[
u_n^{(j)}(x) \coloneqq \partial_j(\calK_A f_n)(x) = \sum_{k=1}^d\partial_kf_n(Ax)\cdot\partial_jA_k(x) = \sum_{k=1}^d(\calK_A\partial_kf_n)(x)\cdot\partial_jA_k(x),\qquad x\in \Omega_X.
\]
For $j\in [1:d]$ consider $u^{(j)}\coloneqq\sum_{k=1}^d \calK_A(\partial_kf)\cdot\partial_jA_k$. Then, Step 1 (applied to $\partial_kf,\partial_kf_n\in L^2(\Omega_\domY)$) and the boundedness of the $\partial_jA_k$ on $\Omega_X$ show that $u_n^{(j)},u^{(j)} \in L^2(\Omega_\domX)$. For fixed $j \in [1:d]$ we obtain
\begin{align*}
\big\|u^{(j)} - u_n^{(j)}\big\|_{L^2(\domX)}
&\le \sum_{k=1}^d\big\|\big[\calK_A(\partial_kf) - \calK_A(\partial_kf_n)\big]\cdot\partial_jA_k\big\|_{L^2(\Omega_X)}\\
&\le \|\det DA^{-1}\|_{C_b(\Omega_Y)}^{1/2}\cdot\sum_{k=1}^d\norm{\partial_jA_k}_{C(\domX)}\norm{\partial_kf - \partial_kf_n}_{L^2(\Omega_\domY)},
\end{align*}
which tends to zero as $n \to \infty$, since $f_n\to f$ in $H^1(\Omega_\domY)$.
In particular, for any test function $\vphi \in C^\infty_c(\Omega_X)$,
\begin{align*}
\int_{\Omega_X} \koopman f\cdot\partial_j\vphi\,dx 
&= \lim_{n\to\infty}\int_{\Omega_X}\koopman f_n\cdot\partial_j\vphi\,dx = -\lim_{n\to\infty}\int_{\Omega_X}u_n^{(j)}\cdot\vphi\,dx = -\int_{\Omega_\domX}u^{(j)}\cdot\vphi\,dx.
\end{align*}
Therefore, $\koopman f\in H^1(\Omega_X)$ with weak partial derivatives $\partial_j\koopman f = u^{(j)}$, $j\in [1:d]$.

\smallskip
\noindent {\bf Step 3.} \cref{eq:Koop_pres_reg} with $\sigma\in\N$, $\sigma\le m$. We prove the claim by induction. For $\sigma\in\{0,1\}$, the claim has already been proven in Steps 1 and 2. Let $\sigma\ge 2$, suppose that $\calK_A : H^{\sigma-1}(\Omega_Y)\to H^{\sigma-1}(\Omega_X)$ is well defined, and consider $f \in H^\sigma(\Omega_\domY)$. By Step 2, we have $\koopman f \in H^1(\Omega_\domX)$ with
\begin{equation}\label{eq:repr-chain-rule}
\partial_j\koopman f = \sum_{k=1}^d \left(\koopman {\partial_kf}\right)\partial_jA_k,\qquad j\in[1:d].
\end{equation}
Since $\partial_kf\in H^{\sigma-1}(\Omega_\domY)$, we obtain $\koopman \partial_kf\in H^{\sigma-1}(\Omega_\domX)$ for $k\in[1:d]$. By assumption, for all $j,k\in[1:d]$, we have $\partial_jA_k\in C^{\sigma-1}(\Omega_\domX)$ with bounded derivatives up to order $\sigma-1$. Hence, we conclude $\partial_j\koopman f\in H^{\sigma-1}(\Omega_\domX)$ for all $j\in[1:d]$ and thus $\koopman f \in H^\sigma(\Omega_\domX)$.

\smallskip
\noindent {\bf Step 4.} Boundedness for $m$. By Step 3, we know that $\calK_A : H^m(\Omega_Y)\to H^m(\Omega_X)$ is well defined. Since $\Omega_X$ and $\Omega_Y$ being Lipschitz domains satisfy the uniform cone condition (cf.\ \cref{app:boundary}) and $m > d/2$, it follows from \cite[Theorems 1 and 2]{Hegland1986} that $H^m(\Omega_X)$ and $H^m(\Omega_Y)$ are native spaces of some kernel functions, respectively. Hence, \cref{lem:koopman-is-continuous-on-native-spaces} and the subsequent discussion imply that $\calK_A : H^m(\Omega_Y)\to H^m(\Omega_X)$ is indeed a bounded operator.

\smallskip
\noindent {\bf Step 5.} The general case. Let $\sigma\in (0,m)$. From the previous steps, we know that $\calK_A : H^m(\Omega_Y)\to H^m(\Omega_X)$ is well defined and bounded. Moreover, the same holds for $\calK_A : L^2(\Omega_Y)\to L^2(\Omega_X)$ (see Step 1). By \cite[Theorem 14.2.3]{Brenner2008}, $H^\sigma(\Omega)$ coincides with the real interpolation space (with equivalent norms\footnote{In fact, $H^\sigma(\R^d) = [L^2(\R^d),H^{m}(\R^d)]_\theta$ with {\em equal} norms, see \cite[Theorem B.7]{McLe00}. However, even for Lip\-schitz domains $\Omega$, the norms of $H^\sigma(\Omega)$ and $[L^2(\Omega),H^{m}(\Omega)]_\theta$ are in general \emph{not} equal, but only equivalent. This has been shown in \cite{ChanHewe15}, proving \cite[Theorem B.8]{McLe00} wrong.}):
\begin{equation}\label{eq:h-sigma-as-interpolation-space}
H^\sigma(\Omega) = \big[L^2(\Omega), H^m(\Omega)\big]_\theta, \qquad\theta = \sigma/m\in (0,1),
\end{equation}
for $\Omega\in\{\Omega_X,\Omega_Y\}$. 
Hence, \cite[Lemma 22.3]{Tartar2007} implies that also
\begin{equation*}
\koopman : \big[L^2(\Omega_\domY), H^{m}(\Omega_\domY)\big]_\theta \to \big[L^2(\Omega_\domX), H^{m}(\Omega_\domX)\big]_\theta
\end{equation*}
is well defined and bounded. Thus, the characterization \cref{eq:h-sigma-as-interpolation-space} yields the result.
\end{proof}

\begin{remark}
In \cref{thm:koopman-faa-di-bruno} in the Appendix we provide explicit bounds on the operator norm of $\calK_A : H^\sigma(\Omega_Y)\to H^\sigma(\Omega_X)$ for $\sigma\in\N_0$.
\end{remark}

Finally, we obtain the desired result for the Koopman operator on the native spaces of the Wendland functions.

\begin{corollary}\label{c:invariance_final}
Let $d\in\N$, $k\in\N_0$, where $d\ge 3$ if $k = 0$. Assume in addition that $A\in C_b^{\ceil{\sigma_{d,k}}}(\Omega_X,\R^d)$. Then the Koopman operator 
\[
\calK_A : \native_{\Phi_{d,k}}(\Omega_\domY)\to\native_{\Phi_{d,k}}(\Omega_\domX)
\]
is well defined and bounded.
\end{corollary}
\begin{proof}
This is an immediate consequence of the combination of \cref{thm:koopman-well-defined}
and the identification in \cref{thm:native-spaces-sobolev-spaces} of the Wendland native spaces $\calN_{\Phi_{d,k}}(\Omega)$ with the Sobolev spaces $H^{\sigma_{d,k}}(\Omega)$, $\Omega \in \{\Omega_\domX,\Omega_\domY\}$, with equivalent norms. Note that $\sigma_{d,k} = (d+1)/2 + k > d/2$.
\end{proof}

\section{Uniform error bounds}\label{sec:uniform-error-bounds}
In this section, we prove our main theorem, containing estimates on the errors $\calK_A - S_\calX\calK_AS_\calY$ and $\calK_A - S_\calX\calK_A$, by combining the results from the previous sections. To this end, we recall the projection operator $S_\calX : \calN(X)\to V_\calX$ from \cref{subsec:regression},
\begin{align*}
S_\calX f = \sum_{i=1}^N\left(\bbK_{\calX}^{-1}f_\calX\right)_i \Phi_{x_i},
\end{align*}
where $f\in\calN(\Omega_X)$ and $f_{\calX i} = f(x_i)$, $i\in[1:N]$. Note that if one considers a Lagrange basis $\{\Phi_{x_i}^L\}_{i=1}^N$ of $V_\calX$, where $\Phi_{x_i}^L(x_j) = \delta_{ij}$ for $i,j\in[1:N]$, then we have
\begin{align*}
S_\calX f = \sum_{i=1}^N f_{\calX i}\Phi_{x_i}^L.
\end{align*}
In what follows, we provide an estimate on $\|I-S_\calX\|_{\calN_{\Phi_{d,k}}(\Omega)\to C_b(\Omega)}$ for the native space of Wendland functions as introduced in \cref{sec:wendland}. Denoting the \emph{fill distance} of a set $\calX = \{x_i\}_{i=1}^N\subset\Omega$ of sample points by
\begin{align*}
h_{\calX} := \sup_{x\in\Omega}\,\min_{1\le i\le N} \|x-x_i\|_2,
\end{align*}
we restate the interpolation error estimate \cite[Theorem 11.17]{Wend04}.

\begin{theorem}\label{thm:intest}
Let $d\in\N$, $k\in\N_0$, and assume that the bounded domain $\Omega\subset\R^d$ satisfies the interior cone condition, cf.\ {\rm\cref{app:boundary}}. Then there are constants $C,h_0>0$ \braces{depending on $d$, $k$, and $\Omega$} such that for every set $\calX = \{x_i\}_{i=1}^N\subset\Omega$ of sample points with $h_{\calX}\le h_0$ and all $\alpha\in\N_0^d$, $|\alpha|\le k$, we have
\begin{align*}
\|D^\alpha f - D^\alpha (S_\calX f)\|_{C_b(\Omega)}\le C h_{\calX}^{k+1/2-|\alpha|} \|f\|_{\calN_{\Phi_{d,k}}(\Omega)},\qquad f\in\calN_{\Phi_{d,k}}(\Omega).
\end{align*}
In particular,
\begin{align}\label{e:I-S}
\|I-S_\calX\|_{\calN_{\Phi_{d,k}}(\Omega)\to C_b(\Omega)}\le C h_{\calX}^{k+1/2}.
\end{align}
\end{theorem}

\noindent We now may state the main result of the paper, leveraging \cref{thm:err_est} in the particular case of native spaces induced by Wendland functions. Therein, the interpolation errors present in the upper bounds of \cref{thm:err_est} vanish if the fill distance tends to zero with a particular rate corresponding to the smoothness of the RBFs.

\begin{theorem}\label{thm:main}
Let $\Omega_X$ and $\Omega_Y$ have a Lipschitz boundary, let $d\in\N$, $k\in\N_0$, where $d\ge 3$ if $k = 0$, and assume that $A\in C_b^{m}(\Omega_X,\R^d)$, where $m = \ceil{\sigma_{d,k}}$. Then there exists a constant $C>0$ such that for any set of centers $\calX = \{x_i\}_{i=1}^N\subset\Omega_X$ we have
\begin{align*}        
\|\calK_A - \wh \calK_A\|_{\calN_{\Phi_{d,k}}(\Omega_Y)\to C_b(\Omega_X)}\,\le\,C h_{\calX}^{k+1/2}.
\end{align*}
Moreover,  there are constants $C_1,C_2 > 0$ such that for any two sets $\calX = \{x_i\}_{i=1}^N\subset\Omega_X$ and $\calY = \{y_j\}_{j=1}\subset\Omega_Y$ of centers,
\begin{align*}        
\|\calK_A - \wh \calK_A^\calY\|_{\calN_{\Phi_{d,k}}(\Omega_Y)\to C_b(\Omega_X)}\,\le\,C_1 h_{\calY}^{k+1/2} + C_2 h_{\calX}^{k+1/2}.
\end{align*}
\end{theorem}
\begin{proof}
By \cref{c:invariance_final}, we have $\calK_A\calN_{\Phi_{d,k}}(\Omega_Y)\subset\calN_{\Phi_{d,k}}(\Omega_X)$. Hence, we may invoke \cref{thm:err_est}, which gives
\begin{align*}
\|\calK_A &- S_\calX\calK_A S_\calY\|_{\calN_{\Phi_{d,k}}(\Omega_Y)\to C_b(\Omega_X)}\\
&\le\|I-S_\calY\|_{\calN_{\Phi_{d,k}}(\Omega_Y)\to C(\Omega_Y)} + \|I-S_\calX\|_{\calN_{\Phi_{d,k}}(\Omega_X)\to C_b(\Omega_X)}\cdot\|\calK_A\|_{V_\calY\to\calN_{\Phi_{d,k}}(\Omega_X)}.
\end{align*}
The terms $\|I-S_\calY\|_{\calN_{\Phi_{d,k}}(\Omega_Y)\to C(\Omega_Y)}$ and $\|I-S_\calX\|_{\calN_{\Phi_{d,k}}(\Omega_X)\to C_b(\Omega_X)}$ may now be bounded using \cref{thm:intest}, and
\[
\|\calK_A\|_{V_\calY\to\calN_{\Phi_{d,k}}(\Omega_X)}\le \|\calK_A\|_{\calN_{\Phi_{d,k}}(\Omega_Y)\to\calN_{\Phi_{d,k}}(\Omega_X)} < \infty
\]
by \cref{c:invariance_final}. Similarly,
\begin{align*}
\|\calK_A - S_\calX\calK_A\|_{\calN_{\Phi_{d,k}}(\Omega_Y)\to C_b(\Omega_X)}
\,\le\,\|I - S_\calX\|_{\calN_{\Phi_{d,k}}(\Omega_X)\to C_b(\Omega_X)}\cdot\|\calK_A\|_{\calN_{\Phi_{d,k}}(\Omega_Y)\to\calN_{\Phi_{d,k}}(\Omega_X)},
\end{align*}
which can be bounded in the same way.
\end{proof}

\begin{remark}
Let us take a closer look at the constants in \cref{thm:main}. 
They are comprised of the constant $C$ from \eqref{e:I-S} and the norms $\|\calK_A\|_{V_\calY\to\calN_{\Phi_{d,k}}(\Omega_X)}$, $\|\calK_A\|_{\calN_{\Phi_{d,k}}(\Omega_Y)\to\calN_{\Phi_{d,k}}(\Omega_X)}$. The constant $C$ stems from \cite[Theorem 11.17]{Wend04} and 
depends on the geometry of the domain $\Omega\in\{\Omega_X,\Omega_Y\}$ and is not specified in \cite{Wend04} either. The operator norms can surely not be estimated using \cref{e:sup} and \cref{e:sup2} which are rather qualitative results. On the other hand, we have \cref{thm:koopman-faa-di-bruno} in the Appendix bounding the norm $\|\calK_A\|_{H^{\sigma_{d,k}}(\Omega_Y)\to H^{\sigma_{d,k}}(\Omega_X)}$, at least for odd $d$. If we knew the equivalence constants for the norms of $H^{\sigma_{d,k}}(\Omega)$ and $\calN_{\Phi_{d,k}}(\Omega)$, $\Omega\in\{\Omega_X,\Omega_Y\}$, we could infer a bound on $\|\calK_A\|_{\calN_{\Phi_{d,k}}(\Omega_Y)\to\calN_{\Phi_{d,k}}(\Omega_X)}$. But such constants are not known. Finding them would involve the determination of (a) the constants in \cite[Theorem 10.35]{Wend04} and (b) the norm of an extension operator $E : H^{\sigma_{d,k}}(\Omega)\to H^{\sigma_{d,k}}(\R^d)$.
\end{remark}

For multi-step predictions we obtain in the same way from \cref{thm:multistep} and the previous \cref{thm:main}:

\begin{corollary}
Let the assumptions of \cref{thm:main} hold. Then there is a constant $\tilde{C}_n > 0$ such that\begin{align*}        
\|\calK_A^n - \wh \calK_A^n\|_{\calN_{\Phi_{d,k}}(\Omega_Y)\to C_b(\Omega_X)}\,\le\,\wt C_n h_{\calX}^{k+1/2}.
\end{align*}
\end{corollary}
Of course, $\wt C_n$ grows with increasing $n$, cf.\ \cref{r:multi_step_might _be_bad}.

\section{Numerical Examples}\label{sec:numerics}
In this section, we illustrate the theoretical results by means of two numerical examples: a Duffing oscillator in two space dimensions and a Lorenz system in three space dimensions.

\subsection{Duffing oscillator}
We consider the dynamics of the Duffing oscillator
\begin{align}\label{eq:duffing}
\begin{pmatrix}
\dot x_1\\
\dot x_2
\end{pmatrix} = \begin{pmatrix}
x_2\\x_1 - 3x_1^3
\end{pmatrix}.
\end{align}
from $\Omega_X = (-2,2)^2$ into $\Omega_Y = (-3,3)^2$. We train our model for the two coordinate maps $f_i(x_1,x_2) = x_i$, $i=1,2$, as observables on uniform grids $\calX = (\delta\,\Z\times\delta\,\Z)\cap [-2,2]^2$ and $\calY = (\delta\,\Z\times\delta\,\Z)\cap [-3,3]^2$ with mesh size $\delta\in\{0.2,0.1,0.05,0.025\}$\footnote{The corresponding fill distances are then given by $h_\calX = h_\calY = \delta/\sqrt 2$.} and perform one-step and multi-step predictions of the Duffing oscillator's flow with time step $\Delta t = 0.02$. 
The validation is performed on a uniform grid in $\Omega_X$ having its points in the centers of the squares indicated by the points of the finest training grid. This ensures that the validation grid does not have any points in common with the interpolation grids $\calX$ on which the error is zero (at least for the approximant $S_\calX\calK_A$).

\smallskip
\noindent\textbf{One-step predictions.} The $L^\infty$- and $L^2$-errors for the two approximations $\wh\calK_A$ and $\wh\calK_A^\calY$ on the respective validation grids are depicted in \cref{tab:err_val_duffing}. We observe different convergence rates corresponding to the smoothness of the kernel. Whereas the kernel with smoothness degree $k=1$ performs better than higher smoothness degrees for $\delta = 0.2$, this changes for smaller fill distances, as predicted by our analysis of Theorem~\ref{thm:main}. Moreover, we observe that asymptotically the two methods behave similarly, as the error values in the respective two right columns of \cref{tab:err_val_duffing} coincide. We note that the numerically observed convergence rates are slightly higher than the proved bounds. E.g.\ for $k=1$ we observe quadratic convergence in Table~\ref{tab:err_val_duffing}, whereas the bound of Theorem~\ref{thm:main} yields the rate $3/2$. This higher order is common in kernel-based methods and may stem from the high regularity of both the flow and the considered coordinate functions.

\begin{figure}[ht!]
\centering
\includegraphics[height = 0.23\textheight]{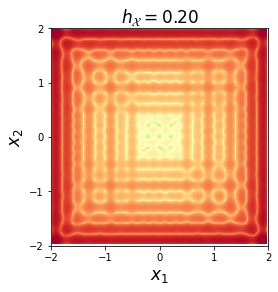}
\includegraphics[height = 0.23\textheight]{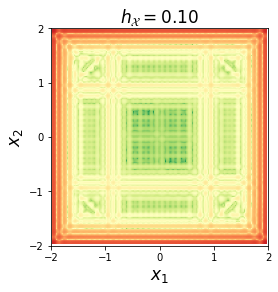}
\includegraphics[height = 0.23\textheight]{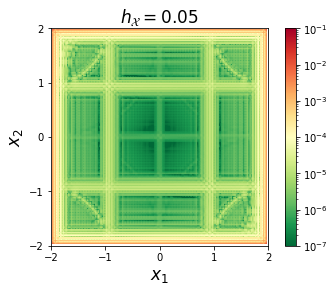}
\caption{$L^\infty$-error for the predicted Duffing oscillator dynamics using the approximant $\wh\calK_A^\calY = S_\calX\calK_A S_\calY$ at smoothness degree $k=1$ with mesh sizes $\delta\in\{0.2,0.1,0.05\}$ from left to right.
}\label{fig:duffing_one}
\end{figure}

\cref{fig:duffing_one} shows intensity plots of the errors
\[
\big\|\big(\calK_A f_1,\calK_A f_2\big)(x_1,x_2) - \big(\wh\calK_A^\calY f_1,\wh\calK_A^\calY f_2\big)(x_1,x_2)\big\|_2,\qquad (x_1,x_2)\in\Gamma_V,
\]
for all considered mesh sizes and smoothness $k=1$. It is readily inferred that the error grows roughly with $\max\{|x_1|,|x_2|\}$ and takes its peak near the boundary of the box. We note that the horizontal and vertical patterns stem from the validation set being a Cartesian product of one-dimensional grids.

\begin{figure}[ht!]\label{fig:Y=X}
\centering
\includegraphics[height = 0.23\textheight]{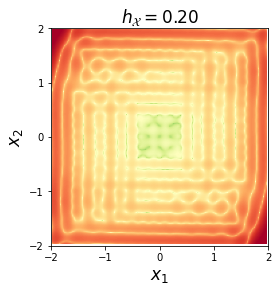}
\includegraphics[height = 0.23\textheight]{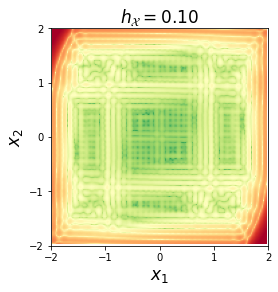}
\includegraphics[height = 0.23\textheight]{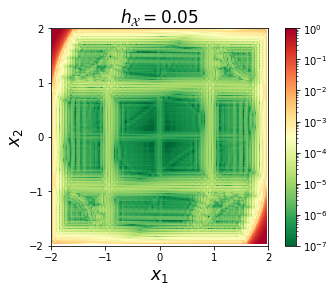}
\label{fig:duffing_one_YgleichX}
\caption{Case $\calY = \calX$: $L^\infty$-error for the predicted Duffing oscillator dynamics using the approximant $S_\calX\calK_A S_\calX$ at smoothness degree $k=1$ with mesh sizes $\delta\in\{0.2,0.1,0.05\}$ from left to right.
}
\end{figure}

In \cref{tab:Y=X} and \cref{fig:Y=X}, we also have depicted the error for the approximation $\wh\calK_A^\calX$ (i.e., $\calY = \calX$). In \cref{fig:Y=X}, we observe that whereas the error in the interior of the domain decreases for smaller fill distances, the dynamics at the corners can not be captured. This leads to a high $L^\infty$-error as reported in \cref{tab:Y=X} and clearly shows that one has to take account of the dynamics in the choice of the $y_i$ if $\Omega_X$ is not invariant under the flow.

\begin{table}[ht!]
\caption{\noindent $L^\infty$-errors \braces{top} and $L^2$-errors \braces{bottom} for the predicted Duffing oscillator dynamics, using the approximants $S_\calX\calK_A$ \braces{left} and $S_\calX\calK_A S_\calY$ \braces{right} with four different fill distances and three different smoothness degrees.} \label{tab:err_val_duffing}
{\footnotesize
\begin{center}
\begin{tabular}[h]{|c|c|c|c|c|}
\hline
$k$ / $\delta$ & 0.2 & 0.1 & 0.05 & 0.025 \\
\hline
1 & 0.15412 & 0.04152 & 0.00982 & 0.00211\\
2 & 0.16564 & 0.02621 & 0.00333 & 0.00033\\
3 & 0.20054 & 0.02154 & 0.00158 & 0.00008\\
\hline
\end{tabular}
\hspace*{.5cm}
\begin{tabular}[h]{|c|c|c|c|c|}
\hline
$k$ / $\delta$ & 0.2 & 0.1 & 0.05 & 0.025 \\
\hline
1 & 0.15283 & 0.04148 & 0.00982 & 0.00211\\
2 & 0.16468 & 0.02620 & 0.00333 & 0.00033\\
3 & 0.19964 & 0.02162 & 0.00158 & 0.00008\\
\hline
\end{tabular}\\[1em]
\begin{tabular}[h]{|c|c|c|c|c|}
\hline
$k$ / $\delta$ & 0.2 & 0.1 & 0.05 & 0.025 \\
\hline
1 & 0.62284 & 0.04950 & 0.00361 & 0.00031\\
2 & 0.70755 & 0.03469 & 0.00142 & 0.00006\\
3 & 0.88888 & 0.03019 & 0.00074 & 0.00001\\
\hline
\end{tabular}
\hspace*{.5cm}
\begin{tabular}[h]{|c|c|c|c|c|}
\hline
$k$ / $\delta$ & 0.2 & 0.1 & 0.05 & 0.025 \\
\hline
1 & 0.61540 & 0.04933 & 0.00361 & 0.00030\\
2 & 0.70539 & 0.03469 & 0.00142 & 0.00006\\
3 & 0.88739 & 0.03020 & 0.00074 & 0.00002\\
\hline
\end{tabular}
\end{center}
}
\end{table}

\begin{table}[ht!]\label{tab:Y=X}
\caption{\noindent Case $\calY = \calX$: $L^\infty$-errors for the predicted Duffing oscillator dynamics using the approximant $S_\calX\calK_A S_\calX$ with four different fill distances and three different smoothness degrees}\label{tab:err_val_duffing_YgleichX}
{\footnotesize
\begin{center}
\begin{tabular}[h]{|c|c|c|c|c|}
\hline
$k$ / $\delta$ & 0.2 & 0.1 & 0.05 & 0.025 \\
\hline
1 & 2.05147 & 1.91387 & 1.82118 & 1.76882\\
2 & 2.29052 & 2.03404 & 1.84288 & 1.72993\\
3 & 2.52979 & 2.21211 & 1.92156 & 1.73675\\
\hline
\end{tabular}
\end{center}
}
\end{table}

\smallskip
\noindent\textbf{Multi-step predictions.} In Table~\ref{tab:multi_rel} we illustrate the relative errors for multi-step predictions in view of the error bound of Theorem~\ref{thm:multistep}. To ensure that the trajectories remain in $\Omega_X=(-2,2)^2$ also for multiple steps, we validate the results for initial states contained in $[-1,1]^2$. We observe that, even for 50 prediction steps, the relative maximal error stays below one percent for the mesh size $\delta = 0.05$ and below five percent for the mesh size $\delta = 0.1$.
\begin{table}[ht!]
    	\caption{\noindent Relative $L^\infty$-errors for $S_\calX\calK_A S_\calX$ with $k=1$.}
    	\footnotesize
    		\begin{center}
    			\begin{tabular}[h]{|c|c|c|c|c|}
    				\hline
    				Steps / $\delta$ & 0.2 & 0.1 & 0.05 \\
    				\hline
    				1 & 0.0038 & 0.00049 & 0.000082\\
    				5 & 0.014 & 0.0034 & 0.00082\\
    				10 & 0.077 & 0.0072 & 0.00095\\
                    15 & 0.14 & 0.016 & 0.0011\\
                    20 & 0.16 & 0.03& 0.0011 \\
                    30 & 0.13& 0.021 & 0.0021 \\
                    40 & 0.176 & 0.034 & 0.0033\\
                    50 & 0.178 & 0.04 & 0.0027\\
    				\hline
    			\end{tabular}
    		\end{center}
\label{tab:multi_rel}
\end{table}

\begin{figure}[htb]
\centering
\includegraphics[height = 0.22\textheight]{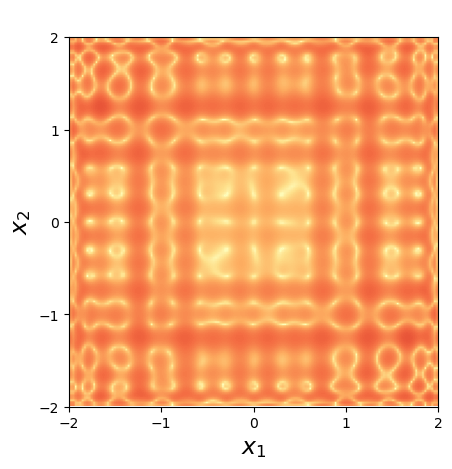}
\includegraphics[height = 0.22\textheight]{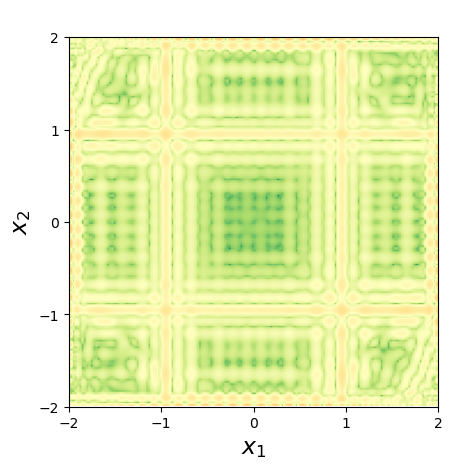}
\includegraphics[height = 0.22\textheight]{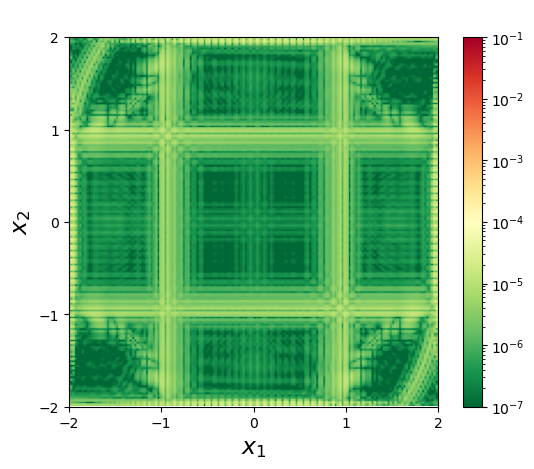}
\caption{$L^\infty$-error for the predicted Duffing oscillator dynamics using the approximant $S_\calX\calK_A S_\calY$ at smoothness degree $k=1$ with meshes of $N\in\{441,1681,6561\}$ Chebyshev nodes, from left to right.
}
\label{fig:cheby}
\end{figure}

\smallskip
\noindent\textbf{Non-uniform meshes.} Figure~\ref{fig:duffing_one} shows that the maximal error occurs at the boundary of the domain. This effect is common in interpolation-based approaches and in our application additionally amplified by the compact support of the radial basis functions. Thus, we briefly investigate grids which become finer towards the boundary. As a simple approach to generate such a mesh in two dimensions we choose a Cartesian product of one-dimensional Chebyshev nodes. In Table~\ref{tab:cheby} we observe that the corresponding maximal error decreases by more than one order of magnitude when compared with the uniform grid, cf.\ Table~\ref{tab:err_val_duffing}. Note that the columns of Table~\ref{tab:cheby} correspond to the respective columns in Table~\ref{tab:err_val_duffing} for $k\in \{0.2,0.1,0.05\}$ in the sense that the associated meshes have the same number of grid points $N$, respectively. Further, in contrast to the relatively high errors at the boundary in Figure~\ref{fig:duffing_one}, the non-uniform meshes lead to a more homogeneous error profile as depicted in Figure~\ref{fig:cheby}.

\begin{table}[ht!]
\caption{\noindent $L^\infty$-errors for the predicted Duffing oscillator dynamics using the approximants $S_\calX\calK_A$ \braces{left} and $S_\calX\calK_A S_\calY$ \braces{right} with a Cartesian product of Chebyshev nodes.} 
\footnotesize
\begin{center}
\begin{tabular}[h]{|c|c|c|c|}
\hline
$k$ / $N$ & 441 & 1681 & 6561 \\
\hline
1 & 0.013 & 0.0023 & 5.1$\cdot 10^{-4}$ \\
2 & 0.0076 & 0.0060 & 3.0$\cdot 10^{-5}$\\
3 & 0.012 & 0.00015 & 4.4$\cdot 10^{-6}$\\
\hline
\end{tabular}
\hspace*{.5cm}
\begin{tabular}[h]{|c|c|c|c|}
\hline
$k$ / $N$ & 441 & 1681 & 6561 \\
\hline
1 & 0.018 & 0.0024 & 5.3$\cdot 10^{-4}$\\
2 & 0.011 & 0.00059 & 2.7$\cdot 10^{-5}$\\
3 & 0.014 & 0.00014 & 5.5$\cdot10^{-6}$\\
\hline
\end{tabular}
\end{center}
\label{tab:cheby}
\end{table}

\subsection{Lorenz system}
In a second numerical experiment, we consider the chaotic Lorenz~63 system, which is given by 
\begin{align}\nonumber 
\dot{x}_1 &= \sigma (x_2-x_1), \\
\dot{x}_2 &= x_1(\rho-x_3)-x_2, \\
\dot{x}_3 &= x_1x_2 - \beta x_3, \nonumber
\end{align}
where $\sigma = 10$, $\rho = 28$, and $\beta = 8/3$. Here, we let $\Omega_X = (0,1)^3$ and fill the domain with four different meshes, each being a Cartesian product of one-dimensional uniform meshes with mesh sizes $\delta\in\{0.2,0.1,0.05,0.025\}$ as before, this time in three dimensions. After simulating the flow with time step $\Delta t = 0.02$, the domain $\Omega_Y$ is chosen as the smallest box including all propagated points; the $\calY$-grids are chosen uniform with fill distances $h_\calY = 2h_\calX$, respectively. The validation grid is again chosen as the union of the centers of the cube cells defined by the $\calX$-grid in order to keep a maximum distance from the interpolation points.

The results are collected in \cref{tab:err_val_lorenz} (top). Contrary to the previous example, where the performance of the two proposed methods was similar, we observe that in particular for large fill distances, the approximation with $S_\calX \calK_A$ leads to smaller errors compared to $S_\calX\calK_A S_\calY$.

In addition, we evaluate our model on the larger domain $\Omega_X = (0,4)^3$ with uniform meshes having the same number of grid points $N$ as the corresponding meshes for the domain $\Omega_X = (0,1)^3$, respectively. In the middle row of \cref{tab:err_val_lorenz}, we provide the corresponding errors. Whereas the error for a small number of grid points is large as depicted in the first two columns, the error significantly decreases when using more sample points and thus reducing the fill distance. In the bottom row of \cref{tab:err_val_lorenz}, we additionally report the errors on the domain $\Omega_X = (0,4)^3$ using a Cartesian product of Chebyshev nodes. We observe that in particular for the highest number of sample points, the error is reduced by more than one order of magnitude compared to the uniform grid in the middle row.

\begin{table}[ht!]
\caption{\noindent $L^\infty$-errors for the predicted dynamics of the Lorenz 63 model using the approximants $S_\calX\calK_A$ \braces{left} and $S_\calX\calK_A S_\calY$ \braces{right}, both with uniform meshes on $[0,1]^3$ \braces{top} and on $[0,4]^3$ with coarser grids \braces{middle} and with a Cartesian product of Chebyshev nodes \braces{bottom}.}\label{tab:err_val_lorenz}
\footnotesize
\begin{center}
\begin{tabular}[h]{|c|c|c|c|c|}
\hline
$k$ / $\delta$ & 0.2 & 0.1 & 0.05 & 0.025 \\
\hline
1 & 0.12825 & 0.02651 & 0.00550 & 0.00118\\
2 & 0.14732 & 0.01737 & 0.00172 & 0.00017\\
3 & 0.20082 & 0.01504 & 0.00085 & 0.00004\\
\hline
\end{tabular}
\hspace*{.5cm}
\begin{tabular}[h]{|c|c|c|c|c|}
\hline
$k$ / $\delta$ & 0.2 & 0.1 & 0.05 & 0.025 \\
\hline
1 & 0.44745 & 0.12192 & 0.02724 & 0.00518\\
2 & 0.54292 & 0.13025 & 0.01678 & 0.00154\\
3 & 0.57911 & 0.16053 & 0.01394 & 0.00075\\
\hline
\end{tabular}

\smallskip 
\begin{tabular}[h]{|c|c|c|c|c|}
\hline
$k$ / $\delta$ & 0.8 & 0.4 & 0.2 & 0.1 \\
\hline
1 & 5.52559 & 1.62021 & 0.47485 & 0.09736\\
2 & 6.93888 & 1.68853 & 0.54650 & 0.06382\\
3 & 7.33920 & 1.23819 & 0.71427 & 0.05702\\
\hline
\end{tabular}
\hspace*{.5cm}
\begin{tabular}[h]{|c|c|c|c|c|}
\hline
$k$ / $\delta$ & 0.8 & 0.4 & 0.2 & 0.1 \\
\hline
1 & 7.23283 & 2.64435 & 1.26929 & 0.36127\\
2 & 7.42979 & 3.76745 & 1.43802 & 0.40461\\
3 & 7.46865 & 4.81280 & 1.49761 & 0.53682\\
\hline
\end{tabular}

\smallskip
\begin{tabular}[h]{|c|c|c|c|c|}
\hline
$k$ / $N$ & 216 & 1331 & 9261 & 68921 \\
\hline
1 & 5.22314 & 0.48521 & 0.01777 & 0.00586\\
2 & 5.85276 & 1.19269 & 0.02473 & 0.00148\\
3 & 6.51056 & 2.09623 & 0.04442 & 0.00040\\
\hline
\end{tabular}
\hspace*{.5cm}
\begin{tabular}[h]{|c|c|c|c|c|}
\hline
$k$ / $N$ & 216 & 1331 & 9261 & 68921 \\
\hline
1 & 6.61044 & 4.92266 & 0.45278 & 0.04314\\
2 & 6.83380 & 5.90063 & 1.02711 & 0.02426\\
3 & 7.12933 & 6.30247 & 1.93234 & 0.03633\\
\hline
\end{tabular}
\end{center}
\end{table}

\section{Conclusions}\label{sec:conclusions}
We have shown in \cref{prop:mat_repr} that the kernel EDMD approximant of the Koopman operator~$\mathcal{K}_{\mathcal{A}}$ --~typically defined as the solution of a linear regression problem~-- may be equivalently expressed as a compression in native spaces by using interpolation operators. 
Based on this novel representation, we derived the first \textit{uniform} finite-data error estimates for kEDMD. This enabled us to prove convergence in the infinite-data limit with convergence rates depending on the smoothness of 
the dynamics in~\cref{thm:main}.
To this end, we have rigorously shown invariance of a rich class of fractional Sobolev spaces under the Koopman operator --~a key property leveraged in the subsequent analysis. 
These 
are generated by Wendland kernels (compactly-supported radial basis functions of minimal degree) and are particularly attractive from a numerical perspective~\cite{Wend04}.


\bigskip
\appendix

\section{Conditions on the boundary of a bounded domain}\label{app:boundary}
In this paper, we consider several conditions on the regularity of the boundary of a bounded domain $\Omega\subset\R^d$. Recall that a \emph{domain} is a non-empty connected open set.

\newpage
\begin{defi}\label{def:boundary}
Let $\Omega\subset\R^d$ be a bounded domain.
\begin{enumerate}
\item[{\rm (a)}] {\rm(\cite[\paragraf 4.9]{Adams2003})} We say that $\Omega$ has a \emph{Lipschitz boundary} if for every $x\in\partial\Omega$ there is a neighborhood $U_x$ of $x$ such that $U_x\cap\partial\Omega$ is the graph of a Lipschitz-continuous function, i.e.\ there exist a Lipschitz-continuous function $\varphi : \R^{d-1}\to\R$ and a rigid motion \braces{i.e., a rotation plus a translation} $T : \R^d\to\R^d$ such that
\[
U_x\cap\Omega = U_x\cap \big\{T(y) : y\in\R^d,\,y_d < \varphi(y_1,\ldots, y_{d-1})\big\}.
\]
\item[{\rm (b)}] {\rm(\cite{Hegland1986})} $\Omega$ satisfies the \emph{uniform cone condition} if there exists a locally finite countable open cover $\{U_j\}$ of the boundary of $\Omega$ and a corresponding sequence $\{C_j\}$ of finite cones, each congruent to some fixed finite cone $C$, such that
\begin{enumerate}
    \item There exists $M>0$ such that every $U_j$ has diameter less then $M$
    \item $\{x\in\Omega : \dist(x,\partial\Omega) < \delta\}\subset\bigcup_jU_j$ for some $\delta>0$
    \item $Q_j := \bigcup_{x\in\Omega\cap U_j}(x+C_j)\subset\Omega$ for every $j$
    \item For some finite $R$, every collection of $R+1$ of the sets $Q_j$ has empty
intersection.
\end{enumerate}

\item[{\rm (c)}] {\rm(\cite[Definition 3.6]{Wend04})} $\Omega$ satisfies the \emph{interior cone condition} if there are an angle $\theta\in (0,\pi/2)$ and a radius $r>0$ such that the following holds: for every $x\in\Omega$ there exists $\xi(x)\in\R^d$, $\|\xi(x)\|_2=1$, such that the cone
\begin{align*}
C\big(x,\xi(x),\theta,r\big) := \{x+\lambda y : y\in \R^d,\,\|y\|=1,\,y^\top\xi(x)\ge\cos\theta,\,\la\in [0,r]\}
\end{align*}
is contained in $\Omega$.
\end{enumerate}
\end{defi}

In \cref{thm:native-spaces-sobolev-spaces}, the domain needs to have a Lipschitz boundary in the sense of \cite[Definition 1.4.4]{Brenner2008}, which, on bounded domains, is easily seen to be implied by the Lipschitz condition above. To obtain the interpolation estimates in \cref{thm:intest}, the \emph{interior cone condition} is needed, which is precisely the \emph{cone condition} from \cite[Paragraph 4.6]{Adams2003} and is implied by the Lipschitz boundary condition (see \cite[Paragraph 4.11]{Adams2003}).

\section{The inclusion property 
and the Koopman operator norm}\label{s:incprop}
The following proposition characterizes the condition \eqref{eq:invariant} in terms of kernel matrices.

\begin{proposition}\label{prop:inv_charac}
The inclusion property~\cref{eq:invariant} holds if and only if
\begin{align}\label{e:sup}
\sup_{\calX}\big\|\bbK_{\calX}^{-1/2}\bbK_{A(\calX)}\bbK_{\calX}^{-1/2}\big\|_{2\to 2}^{1/2}\,<\,\infty,
\end{align}
where $\calX$ runs through all finite sets $\calX\subset X$, and $\|\cdot\|_{2\to 2}$ denotes the matrix spectral norm. In this case, the operator norm of $\calK_A|_{\calN(Y)} : \calN(Y)\to\calN(X)$ is given by the supremum in \eqref{e:sup}.
\end{proposition}
\begin{proof}
Assume that \eqref{eq:invariant} holds. Then the operator $\calK_A : \calN(Y)\to\calN(X)$ is bounded by \cref{lem:koopman-is-continuous-on-native-spaces}. Hence, we may consider its Hilbert space adjoint $\calK_A^*:\calN(X)\to \calN(Y)$, which by \eqref{eq:koopmanadjoint} satisfies  $\calK_A^*\Phi_x = \Psi_{A(x)}$ for $x\in X$. Hence, as $\linspan\{\Phi_x : x\in X\}$ is dense in $\calN(X)$,
\begin{align*}
\|\calK_A^*\|_{\calN(X)\to\calN(Y)}^2
&= \sup_\calX\,\sup\left\{\Big\|\calK_A^*\sum_{i=1}^n\alpha_i\Phi_{x_i}\Big\|_{\calN(Y)}^2 : n = |\calX|,\,\alpha\in\R^n,\,\Big\|\sum_{i=1}^n\alpha_i\Phi_{x_i}\Big\|_{\calN(X)}^2 = 1\right\}\\
&= \sup_\calX\,\sup\Big\{\alpha^\top\bbK_{A(\calX)}\,\alpha : n = |\calX|,\,\alpha\in\R^n,\,\alpha^\top\bbK_{\calX}\,\alpha = 1\Big\}\\
&= \sup_\calX\,\sup\left\{\alpha^\top\bbK_{A(\calX)}\,\alpha : n = |\calX|,\,\alpha\in\R^n,\,\|\bbK_{\calX}^{1/2}\,\alpha\|_2 = 1\right\}\\
&= \sup_\calX\,\sup\left\{v^\top\bbK_{\calX}^{-1/2}\bbK_{A(\calX)}\,\bbK_{\calX}^{-1/2}v : n = |\calX|,\,v\in\R^n,\,\|v\|_2 = 1\right\}\\
&= \sup_\calX\big\|\bbK_{\calX}^{-1/2}\bbK_{A(\calX)}\,\bbK_{\calX}^{-1/2}\big\|_{2\to 2}.
\end{align*}
The claim now follows from the equality of the norms of $\calK_A$ and $\calK_A^*$.

Conversely, assume that \eqref{e:sup} is satisfied, and let $S$ be the supremum in \eqref{e:sup}. Define the normed spaces
\[
\calN_0(X) = \big(\linspan\{\Phi_x : x\in X\},\|\cdot\|_{\calN(X)}\big)
\quad\text{and}\quad
\calN_0(Y) = \big(\linspan\{\Psi_y : y\in Y\},\|\cdot\|_{\calN(Y)}\big),
\]
and a linear operator $T_0 : \calN_0(X)\to\calN_0(Y)$ by $T_0\Phi_x = \Psi_{A(x)}$ for $x\in X$. Then the calculations above show that $\sup\{\|T_0f\|_{\calN(Y)} : f\in\calN_0(X),\,\|f\|_{\calN(X)}=1\} = S^{1/2}$. Hence, $T_0$ is bounded with norm $S^{1/2}$ and thus admits a bounded extension $T : \calN(X)\to\calN(Y)$. For any $f\in\calN(Y)$ and $x\in X$ we have
\[
(T^*f)(x) = \<T^*f,\Phi_x\>_{\calN(X)} = \<f,T_0\Phi_x\>_{\calN(Y)} = \<f,\Psi_{A(x)}\>_{\calN(Y)} = f(A(x)) = (\calK_A f)(x),
\]
thus $\calK_A f = T^*f\in\calN(X)$, which proves \eqref{eq:invariant}.
\end{proof}

We point out that Proposition \ref{prop:inv_charac} also holds for discontinuous 
kernels $k_X$ and $k_Y$ and even without the strict positive definiteness of~$k_Y$. 
However, we remark the following 
on the strict definiteness of the kernel $k_X$.

\begin{remark}
For non-strictly positive definite kernels $k_X$ (e.g., polynomial kernels), the inverse of $\bbK_\calX$ in \eqref{e:sup} is not defined for some sets $\calX$. Even more, for \eqref{eq:invariant} to hold, it is necessary that
\[
\ker\bbK_\calX\subset\ker\bbK_{A(\calX)}
\]
for all finite sets $\calX\subset X$. To see this, we follow the calculation for $\|\calK_A^*\|_{\calN(X)\to\calN(Y)}^2$ in the proof of \cref{prop:inv_charac} until the third line. Now, suppose that $w\in\ker\bbK_\calX$ such that $\bbK_{A(\calX)}w\neq 0$, and let $v\in\R^n$ such that $\|\bbK_\calX^{1/2}v\|_2=1$. Set $\alpha_\la=\la w+v$, $\la>0$. Then $\|\bbK_\calX^{1/2}\alpha_\la\|_2=1$ for every $\la > 0$, and, hence,
\[
\|\calK_A^*\|_{\calN(X)\to\calN(Y)}\ge \alpha_\la^\top\bbK_{A(\calX)}\alpha_\la = \la^2\<\bbK_{A(\calX)}w,w\> + 2\la\<\bbK_{A(\calX)}w,v\> + \<\bbK_{A(\calX)}v,v\>.
\]
But the right-hand side tends to $\infty$ as $\la\to\infty$, which is a contradiction.

To conclude this remark, we mention that \eqref{eq:invariant} is satisfied in general if and only if $\ker\bbK_\calX\subset\ker\bbK_{A(\calX)}$ holds for all finite sets $\calX$ and $\sup_{\calX}\|(\bbK_\calX^\dagger)^{1/2}\bbK_{A(\calX)}(\bbK_\calX^\dagger)^{1/2}\|_{2\to 2}^{1/2} < \infty$.
\end{remark}

The next corollary shows that the operator norm of $\calK_A : \calN(Y)\to\calN(X)$ is typically larger than one.

\begin{corollary}\label{c:norm_large}
If \eqref{eq:invariant} holds and $c = k_X(x,x) = k_Y(y,y)$ for all $x\in X$ and $y\in Y$, then
\[
\|\calK_A\|_{\calN(Y)\to\calN(X)}^2\ge 1 + \sup_{x,z\in X}\,\frac{|d(x,z)|}{c + \sgn(d(x,z))k_X(x,z)},
\]
where $d(x,z) = k_Y(A(x),A(z)) - k_X(x,z)$. In particular, we have $\|\calK_A\|_{\calN(Y)\to\calN(X)} > 1$, except in the case $k_X = k_Y(A(\cdot),A(\cdot))$, where $\|\calK_A\|_{\calN(Y)\to\calN(X)} = 1$.
\end{corollary}
\begin{proof}
Let $x,z\in X$ be arbitrary, $z\neq x$, and let $\calX = \{x,z\}$. Then, setting $a = k_X(x,z)$ and $w_\pm = \sqrt{c\pm a}$, we have
\[
\bbK_\calX^{-1/2} = \bmat caac^{-1/2} = \frac 1{2w_+w_-}\bmat{w_-+w_+}{w_--w_+}{w_--w_+}{w_-+w_+},
\]
and hence, with $b = k_Y(A(x),A(z))$,
\[
\bbK_{\calX}^{-1/2}\bbK_{A(\calX)}\,\bbK_{\calX}^{-1/2} = \frac{1}{c^2-a^2}\bmat{c^2-ab}{c(b-a)}{c(b-a)}{c^2-ab}.
\]
Setting $\sigma = \sgn(b-a) = \sgn(d(x,z))$, the spectral norm of this matrix equals \begin{align*}
\frac{c^2-ab+c|b-a|}{c^2-a^2}
&= 1 + \frac{a^2-ab+c\sigma(b-a)}{c^2-a^2} = 1 + \frac{\sigma a(\sigma a-c) - \sigma b(\sigma a - c)}{(c-\sigma a)(c + \sigma a)} = 1 + \frac{|b-a|}{c+\sigma a}.
\end{align*}
The claim now follows from \cref{prop:inv_charac}.
\end{proof}

\section{Bounds on the Koopman operator norm}\label{app:constants}
In this section, we prove explicit bounds on the operator norms $\|\calK_A\|_{H^\sigma(\Omega_Y)\to H^\sigma(\Omega_X)}$ for $\sigma\in\N$. We denote by $\mathcal L^{k}(\R^d,\R^n)$ the linear space of all $k$-multilinear mappings $\Psi : (\R^d)^k\to \R^n$. A multilinear map $\Psi\in\calL^{k}(\R^d,\R^n)$ is called {\em symmetric} if $\Psi(Pv) = \Psi(v)$ for any permutation matrix $P\in\R^{d\times d}$ and all $v\in(\R^d)^k$. By $\calL_s^{k}(\R^d,\R^n)$ we denote the set of all symmetric $k$-multilinear maps. 

For the $k$-th total derivative $\calD^kA$ of $A : \Omega_X\to\Omega_Y$ we have $\calD^kA(x)\in\calL^k_s(\R^d,\R^d)$, $x\in\Omega_X$, by setting
\[
\calD^kA(x)(e_{i_1},\ldots,e_{i_k}) = \partial_{i_1}\cdots\partial_{i_k}A(x),
\]
where $e_i$ denotes the $i$-th standard basis vector. For scalar-valued $f : \Omega_Y\to\R$ we have accordingly $f^{(k)}(y) := \calD^k f(y)\in\calL^k_s(\R^d,\R)$. For example, the second derivative of $f$ can be written as $f^{(2)}(y)(x_1,x_2) = x_1^\top H_f(y)x_2$, where $H_f(y)$ denotes the Hessian of $f$ at $y$. 

For $\sigma \in \N$ we denote the set of all partitions $\pi$ of $\{1, \dots, \sigma\}$ by $\Pi_\sigma$. For a partition $\pi \in \Pi_\sigma$, $\abs \pi$ denotes the number of blocks $B$ in the partition $\pi$.

For $\alpha\in\N_0^d$, we denote by $\calP(\alpha)$ the set of all vectors $\beta\in [1:d]^{|\alpha|}$ in which any index $k\in [1:d]$ appears exactly $\alpha_k$ times. For example, for $\alpha = (2,1)$ we have $\calP(\alpha) = \{(1,1,2), (1,2,1), (2,1,1)\}$.

\begin{theorem}\label{thm:koopman-faa-di-bruno}
Let $\sigma\in\N$, and assume in addition that $A\in C_b^\sigma(\Omega_X,\R^d)$. Then we have
\[
\norm{\koopman}_{H^{\sigma}(\domY) \to H^{\sigma}(\domX)}\le
\bigg(\max\bigg\{c_0,\sum_{1\le|\alpha|\le\sigma}S_{\alpha}(A)\bigg\}\bigg)^{1/2},
\]
where $c_0 := \sup_{x\in\Omega_X}|\det DA(x)|^{-1}$, and for $\alpha\in\N_0^d$, $|\alpha|=m$, and any $\beta\in\calP(\alpha)$,
\[
S_\alpha(A) := \sup_{x\in\Omega_X}|\det DA(x)|^{-1}\sum_{k=1}^{m}\sum_{|\alpha'|=k}\Bigg|\sum_{\substack{\pi\in\Pi_m\\|\pi|=k}}\sum_{\beta'\in\calP(\alpha')}\partial_{\beta,B_1}A_{\beta'_1}(x)\cdots\partial_{\beta,B_k}A_{\beta'_k}(x)\Bigg|^2.
\]
Here, $B_i$ denotes the $i$-th block of the partition $\pi$, and $\partial_{\beta,B}$ stands for the operator $\partial_{\beta_{j_1}}\cdots\partial_{\beta_{j_\ell}}$, where $B = \{j_1,\ldots,j_\ell\}$.
\end{theorem}
\begin{proof}
We prove by induction over $\sigma\in\N$ that for $f\in H^\sigma(\Omega_Y)$ we have
\begin{align}\label{eq:claim}
\|\calK_Af\|_{H^\sigma(\Omega_X)}^2\,\le\,c_0\|f\|_{L^2(\Omega_Y)}^2 + \sum_{k=1}^\sigma\bigg(\sum_{k\le|\alpha|\le\sigma}S_\alpha(A)\bigg)\sum_{|\gamma|=k}\|D^\gamma f\|_{L^2(\Omega_Y)}^2.
\end{align}
Then the statement of the theorem follows immediately. The anchor for \cref{eq:claim} has been set in the first step of the proof of \cref{thm:koopman-well-defined}, where it was shown that $\|\calK_A\|_{L^2(\Omega_Y)\to L^2(\Omega_X)}^2\le c_0$. Next, let $\sigma\in\N$, $\sigma\ge 1$. For any $\alpha\in\N_0^d$ with $|\alpha|=\sigma$, we have $D^\alpha = \partial_\beta$ with some (any) $\beta\in\calP(\alpha)$. Now, consider Fa\`a di Bruno's formula in combinatorial form (see \cite[p. 219]{Johnson02}):
\begin{equation*}
(f \circ A)^{(\sigma)}(x)(z_1,\ldots,z_\sigma) = \sum_{\pi \in \Pi_\sigma} f^{({\abs{\pi}})}(A(x))\big(\calD^{\abs{B_1}}A(x)(z_{B_1}), \dots, \calD^{\abs{B_{|\pi|}}}A(x)(z_{B_{|\pi|}}) \big),
\end{equation*}
where $z_{B_j}\in (\R^d)^{|B_j|}$ are the components of $(z_1,\ldots,z_\sigma)$ with indices in the block $B_j$. If we denote $e^\beta = (e_{\beta_1},\ldots,e_{\beta_\sigma})$, then
\begin{align*}
D^\alpha(f\circ A)(x)
&= \partial_\beta(f\circ A)(x) = (f\circ A)^{(\sigma)}(x)(e^\beta)\\
&= \sum_{k=1}^\sigma\sum_{\substack{\pi\in\Pi_\sigma\\|\pi|=k}} f^{(k)}(A(x))\big(\calD^{\abs{B_1}}A(x)(e^\beta_{B_1}), \dots, \calD^{\abs{B_{k}}}A(x)(e^\beta_{B_k})\big)\\
&= \sum_{k=1}^\sigma\sum_{\substack{\pi\in\Pi_\sigma\\|\pi|=k}} f^{(k)}(A(x))\big(\partial_{\beta,B_1}A(x),\dots,\partial_{\beta,B_k}A(x)\big)\\
&= \sum_{k=1}^\sigma\sum_{\substack{\pi\in\Pi_\sigma\\|\pi|=k}}\sum_{\beta'\in [1:d]^k}\partial_{\beta,B_1}A_{\beta'_1}(x)\cdots\partial_{\beta,B_k}A_{\beta'_k}(x)\cdot f^{(k)}(A(x))(e^{\beta'})\\
&= \sum_{k=1}^\sigma\sum_{\substack{\pi\in\Pi_\sigma\\|\pi|=k}}\sum_{\beta'\in [1:d]^k}\partial_{\beta,B_1}A_{\beta'_1}(x)\cdots\partial_{\beta,B_k}A_{\beta'_k}(x)\cdot \partial_{\beta'}f(A(x))\\
&= \sum_{k=1}^\sigma\sum_{|\alpha'|=k}D^{\alpha'}f(A(x))\underbrace{\sum_{\substack{\pi\in\Pi_\sigma\\|\pi|=k}}\sum_{\beta'\in\calP(\alpha')}\partial_{\beta,B_1}A_{\beta'_1}(x)\cdots\partial_{\beta,B_k}A_{\beta'_k}(x)}_{=: A_{\alpha,k,\alpha'}(x)}.
\end{align*}

Thus, we have $D^\alpha(f\circ A) = \sum_{k=1}^\sigma\sum_{|\alpha'|=k}[(D^{\alpha'}f)\circ A]\cdot A_{\alpha,k,\alpha'}$, and we estimate
\begin{align*}
\|D^\alpha(f\circ A)\|_{L^2(\Omega_X)}^2
&= \int_{\Omega_X}\bigg(\sum_{k,\alpha'}[(D^{\alpha'}f)\circ A]\cdot A_{\alpha,k,\alpha'}\bigg)^2\,dx\\
&= \int_{\Omega_Y}\bigg(\sum_{k,\alpha'}(D^{\alpha'}f)\cdot [A_{\alpha,k,\alpha'}\circ A^{-1}]\bigg)^2|\det DA^{-1}|\,dy\\
&= \int_{\Omega_Y}\bigg(\sum_{k,\alpha'}(D^{\alpha'}f)\cdot\underbrace{[A_{\alpha,k,\alpha'}\circ A^{-1}]|\det DA^{-1}|^{1/2}}_{=: \wt A_{\alpha,k,\alpha'}}\bigg)^2\,dy\\
&\le \int_{\Omega_Y}\bigg(\sum_{k,\alpha'}|D^{\alpha'}f|^2\bigg)\bigg(\sum_{k,\alpha'}|\wt A_\alpha,{k,\alpha'}|^2\bigg)\,dy\\
&\le \bigg(\sup_{y\in\Omega_Y}\sum_{k,\alpha'}|\wt A_{\alpha,k,\alpha'}(y)|^2\bigg)\bigg(\sum_{k,\alpha'}\int_{\Omega_Y}|D^{\alpha'}f|^2\,dy\bigg)\\
&\le \bigg(\sup_{x\in\Omega_X}\sum_{k,\alpha'}\frac{|A_{\alpha,k,\alpha'}(x)|^2}{|\det DA(x)|}\bigg)\sum_{1\le|\gamma|\le\sigma}\|D^\gamma f\|_{L^2(\Omega_Y)}^2\\
&= S_\alpha(A)\sum_{1\le|\gamma|\le\sigma}\|D^\gamma f\|_{L^2(\Omega_Y)}^2.
\end{align*}
Therefore, we obtain by assumption that
\begin{align*}
\|\calK_Af\|_{H^\sigma(\Omega_X)}^2
&= \|\calK_Af\|_{H^{\sigma-1}(\Omega_X)}^2 + \sum_{|\alpha|=\sigma}\|D^\alpha(f\circ A)\|_{L^2(\Omega_X)}^2\\
&\le c_0\|f\|_{L^2(\Omega_Y)}^2 + \sum_{k=1}^{\sigma-1}\bigg(\sum_{k\le|\alpha|\le\sigma-1}S_\alpha(A)\bigg)\sum_{|\gamma|=k}\|D^\gamma f\|_{L^2(\Omega_Y)}^2\\
&\hspace*{2.47cm}+ \sum_{|\alpha|=\sigma}S_\alpha(A)\sum_{1\le|\gamma|\le\sigma}\|D^\gamma f\|_{L^2(\Omega_Y)}^2\\
&= c_0\|f\|_{L^2(\Omega_Y)}^2 + \sum_{k=1}^{\sigma}\bigg(\sum_{k\le|\alpha|\le\sigma}S_\alpha(A)\bigg)\sum_{|\gamma|=k}\|D^\gamma f\|_{L^2(\Omega_Y)}^2,
\end{align*}
as we wished to prove.
\end{proof}

In the particular case $\sigma=1$, \cref{thm:koopman-faa-di-bruno} yields the following corollary.

\begin{corollary}\label{cor:koopman-norm-h1}
We have
\[
\norm{\koopman}_{H^1(\Omega_Y)\to H^1(\Omega_X)}\le\max\bigg\{c_0,\sup_{x\in\Omega_X}\frac{\|DA(x)\|_F^2}{|\det DA(x)|}\bigg\}^{1/2}.
\]
\end{corollary}

\begin{remark}
In fact, it is not hard to see that the Frobenius norm  $\|DA(x)\|_F$ above can be replaced by the spectral norm $\|DA(x)\|_{2\to 2}$ (cf.\ \cite{GoldGuro95}).
\end{remark}
\end{document}